\documentclass[journal]{IEEEtran}
\usepackage{cite}
\usepackage{url}
\usepackage{empheq}
\usepackage{float}
\usepackage{dsfont}
\usepackage[hidelinks]{hyperref}

\usepackage{algorithm}
\usepackage{algorithmic}

\usepackage{amsmath,amssymb,amsfonts}
\allowdisplaybreaks[4]
\usepackage{graphicx}
\usepackage{textcomp}
\usepackage{amsmath}
\usepackage{enumerate}
\usepackage{epstopdf}
\usepackage{array}
\usepackage{booktabs}
\usepackage{subfigure}
\usepackage{multirow}
\usepackage{soul, color}
\usepackage[usenames,dvipsnames]{xcolor}
\usepackage[latin1]{inputenc}
\usepackage{cases}

\newtheorem{proposition}{Proposition}

\soulregister\cite7 
\soulregister\citep7 
\soulregister\citet7 
\soulregister\ref7 
\soulregister\pageref7

\usepackage{bbding}

\usepackage{nomencl}
\makenomenclature
\usepackage{etoolbox}
\renewcommand\nomgroup[1]{%
  \item[\bfseries
  \ifstrequal{#1}{A}{Acronyms}{%
  \ifstrequal{#1}{S}{Symbols}{%
  \ifstrequal{#1}{U}{Units}{}}}%
]}

\begin{document}

\title{Real-time Feedback Based Online Aggregate EV Power Flexibility Characterization}

\author{
Dongxiang~Yan,
Shihan~Huang,
and Yue~Chen,~\IEEEmembership{Member,~IEEE}
\thanks{This work was supported by the National Natural Science Foundation of China under Grant No. 52307144 and the Shun Hing Institute of Advanced Engineering, the Chinese University of Hong Kong, through Project RNE-p2-23 (Grant No. 8115071). (Corresponding to Y. Chen)}
\thanks{D. Yan, S. Huang, and Y. Chen are with the Department of Mechanical and Automation Engineering, the Chinese University of Hong Kong, Hong Kong SAR, China (e-mail: dongxiangyan@cuhk.edu.hk, shhuang@link.cuhk.edu.hk, yuechen@mae.cuhk.edu.hk).
}}
\markboth{Journal of \LaTeX\ Class Files,~Vol.~XX, No.~X, Feb.~2019}%
{Shell \MakeLowercase{\textit{et al.}}: Bare Demo of IEEEtran.cls for IEEE Journals}

\maketitle

\begin{abstract}
As an essential measure to combat global warming, electric vehicles (EVs) have witnessed rapid growth. Flexible EVs can enhance power systems' ability to handle renewable generation uncertainties. How EV flexibility can be utilized in power grid operation has captured great attention. However, the direct control of individual EVs is challenging due to their small capacity and large number. Hence, it is the aggregator that interacts with the grid on behalf of the EVs by characterizing their aggregate flexibility. In this paper, we focus on the aggregate EV power flexibility characterization problem. First, an offline model is built to obtain the lower and upper bounds of the aggregate EV power flexibility region. It ensures that any trajectory within the region is feasible. Then, considering that parameters such as real-time electricity prices and EV arrival/departure times are not known in advance, an online algorithm is developed based on Lyapunov optimization techniques. We provide a theoretical bound for the maximum charging delay under the proposed online algorithm.
Furthermore, real-time feedback is designed and integrated into the proposed online algorithm to better unlock EV power flexibility. Comprehensive performance comparisons are carried out to demonstrate the advantages of the proposed method.
\end{abstract}

\begin{IEEEkeywords}
Aggregate flexibility, charging station, electric vehicle, Lyapunov optimization, online algorithm.
\end{IEEEkeywords}

\IEEEpeerreviewmaketitle

\section*{Nomenclature}
\addcontentsline{toc}{section}{Nomenclature}
\begin{IEEEdescription}[\IEEEusemathlabelsep\IEEEsetlabelwidth{$e_{v}^{min},e_{v}^{max}$ }]
\item[$\mathcal{T},t$]{Set of time slots and the index.}
\item[$\mathcal{V},v$]{Set of EVs and the index.}
\item[$\mathcal{G},g$]{Set of EV groups and the index.}
\item[$\mathcal{V}_g$]{Set of EVs in group $g$.}

\item[$t_{v}^a,t_{v}^d$]{EV $v$'s arrival time and departure time.}
\item[$e_{v}^{a},e_{v}^{d}$]{EV $v$'s energy level at time $t_{v}^a,t_{v}^d$.}
\item[$R_{v}$]{EV $v$'s maximum allowed charging delay.}
\item[$p_{v,t}^c$]{EV $v$'s charging power at time $t$.}
\item[$\hat{p}_{d,t}/\check{p}_{d,t}$]{Upper/Lower bound of the aggregate EV power flexibility region.}
\item[$F_t$]{Aggregate power flexibility value at $t$.}
\item[$\pi_t$]{Electricity price at time $t$.}
\item[$p_{v}^{max}$]{EV $v$'s maximum charging power.}
\item[$\Delta t$]{Time interval.}
\item[$\hat{e}_{v,t}/\check{e}_{v,t}$]{Upper/Lower bound of EV $v$'s energy at $t$.}
\item[$e_{v}^{min}/e_{v}^{max}$]{Minimum/Maximum EV energy level.}
\item[$R_g$]{Charging delay of group $g$.}
\item[${Q}_{g,t}$]{Charging task queue backlog of group $g$ at time $t$, and $\hat{Q}_{g,t}/\check{Q}_{g,t}$ for the upper/lower bound of the aggregate flexibility region.}
\item[$x_{g,t}$]{Charging power for EVs of group $g$ at $t$.}
\item[$\hat{x}_{g,t}/\check{x}_{g,t}$]{Upper/Lower bound of $x_{g,t}$.}
\item[$a_{g,t}$]{Arrival rate of EV charging tasks of group $g$ at time $t$.}
\item[$\hat{a}_{g,t}/\check{a}_{g,t}$]{Upper/Lower bound of $a_{g,t}$.}
\item[$a_{v,t}$]{Charging demand of EV $v$ at $t$.}
\item[$\hat{a}_{v,t}/\check{a}_{v,t}$]{Upper/Lower bound of $a_{v,t}$.}
\item[$\hat{t}_{v}^{min}/\check{t}_{v}^{min}$]{Upper/Lower bound of minimum required charging time.}
\item[$\eta_c$]{EV charging efficiency.}
\item[$e_{v}^{cha}$]{EV $v$'s required charging demand.}
\item[$\hat{e}_{v}^{cha}/\check{e}_{v}^{cha}$]{Upper/Lower bound of $e_{v}^{cha}$.}
\item[$\eta_g$]{Adjustment parameter.}
\item[$\mathbb{I}_{\hat{Q}_{g,t}},\mathbb{I}_{\check{Q}_{g,t}}$]{Indicator functions of $\hat{Q}_{g,t}$ and $\check{Q}_{g,t}$.}
\item[$Z_{g,t}$]{Delay-aware virtual queue backlog of group $g$ at time $t$, and $\hat{Z}_{g,t}/\check{Z}_{g,t}$ for the upper/lower bound of the aggregate flexibility region.}
\item[$\delta_{g,max}$]{Maximum charging delay of group $g$.}
\item[$\Theta(t)$] {Concatenated vector of virtual queues.}
\item[$V$]{A parameter in the objective function.}
\item[$\alpha_t$]{Dispatch ratio at time $t$.}
\item[$p_{g,t}^{disp}$]{Dispatch power for group $g$ at time $t$.}
\end{IEEEdescription}

\section{Introduction}
\IEEEPARstart {T}{hanks} to the low carbon emissions, electric vehicles (EVs) have been considered a promising solution to climate change and have proliferated in recent years \cite{wang2018electrical}. However, the uncontrolled charging of a large number of EVs can cause voltage deviation, line overload, and huge transmission losses \cite{munshi2018extract}, threatening the reliability of power systems.
Unlike inelastic loads, the charging power and charging period of EVs are more flexible \cite{gan2012optimal}. Therefore, unlocking the power flexibility of EVs is an important way to reduce the adverse impact of high charging demand on the power grids.

Extensive literature has focused on designing coordinated charging strategies to optimally schedule EV charging.
For example, to promote the use of local renewable generation, a dynamic charging strategy was proposed to allow the EV charging power to dynamically track the PV generation \cite{mouli2016system} and wind generation \cite{yang2018distributed}. To reduce the electricity costs, a deterministic optimal charging strategy was proposed for a home energy management system based on the time-of-use tariffs \cite{wu2017optimal}. 
To address the uncertainties associated with EV charging, reference \cite{wu2018stochastic} proposed a stochastic charging strategy based on a probabilistic model related to EV daily travels.
A multi-stage energy management strategy was developed for a charging station integrated with PV generation and energy storage \cite{yan2018optimized}.
In addition, a pricing mechanism was suggested in \cite{zhao2018real} to guide the EVs for economical charging. A double-layer optimization model was built to reduce the voltage deviation caused by EV charging \cite{dong2018charging}.
While great efforts have been made in determining the optimal EV charging power, it is challenging to directly control a large number of EVs due to the high computational complexity.

To address this problem, some studies characterized the aggregate EV power flexibility. 
Reference \cite{xu2016hierarchical} proposed to derive the aggregate EV charging flexibility region by directly summing up the lower and upper bounds of power and cumulative energy. 
This aggregate EV model was adopted by \cite{zhang2016evaluation} to evaluate the achievable vehicle-to-grid capacity of an EV fleet and by \cite{wang2022evaluation} to quantify the value of EV flexibility in maintaining distribution system reliability. Reference \cite{shi2021optimal} further considered the spatio-temporal distribution of EVs. An EV dispatchable region was proposed to facilitate the participation of charging stations in electricity markets \cite{zhou2021forming}.
In addition to EVs, the aggregate flexibility of thermostatically controllable loads (TCLs)\cite{chen2021scheduling,zhao2017geometric}, distributed energy resources \cite{yi2021aggregate,chen2020aggregate}, and virtual power plants \cite{wang2021aggregate} was also studied.


The above works have provided fruitful techniques for evaluating aggregate EV flexibility in an offline manner, which assumes that the aggregator has complete information about future uncertainty realizations. However, in practice, the future EV arrival/departure time and electricity prices can hardly be obtained accurately in advance. Therefore, an online algorithm based solely on the up-to-date information is necessary. A straightforward approach is the greedy algorithm that decomposes the offline problem into subproblems for each time slot by neglecting the time-coupling constraints \cite{guo2021asynchronous}. However, the result could be far from the offline optimum.
An online model predictive control (MPC) algorithm was proposed to minimize the operation cost of EV charging stations \cite{2019zheng} based on short-term forecasts. A combined robust and stochastic MPC method was developed in \cite{jiao2022online} to handle the uncertain EV charging behaviors and renewable generations.
However, the MPC-based methods still rely on forecasts and can be computationally demanding due to the rolling optimization process.
Unlike the MPC-based methods, the reinforcement learning (RL)-based online charging methods enable charging decisions based on the current observations. Reference \cite{dimitrov2014RL} used a basic Q-learning method to schedule individual EV charging. Linear \cite{wang2019RL} and nonlinear \cite{chics2016RL} function approximators for the Q-table were designed. A deep Q-network (DQN) was used to approximate the optimal action-value function to determine the EV charging strategy \cite{wan2019model}. 
In addition, various deep RL (DRL) algorithms were applied, such as the deep deterministic policy gradient (DDPG) \cite{qiu2020deep} and the proximal policy optimization (PPO) \cite{jin2023optimal}. Reference \cite{zhang2021cddpg} combined DDPG and long short-term memory (LSTM) network to learn the optimal EV sequential charging strategy. Reference \cite{yan2021deep} adopted the soft actor-critic (SAC) DRL algorithm and incorporated the driver behaviors. However, the (deep) RL methods face challenges such as their sensitivity to hyperparameters and limited interpretability.

An alternative approach to developing an online algorithm is based on Lyapunov optimization. Compared to the MPC-based and RL methods, the Lyapunov optimization-based methods are prediction-free and have theoretical performance guarantees \cite{fan2020online}. 
Lyapunov optimization has been used in microgrid control \cite{2017RTMG}, energy storage sharing \cite{zhong2019online}, and data center energy management \cite{yu2018distributed}.
A charging strategy based on Lyapunov optimization was proposed to minimize the total electricity cost \cite{shen2021real}. To ensure the satisfaction of charging requirements, virtual delay queues were introduced \cite{jin2014optimized,zhou2018optimal}.

Though the optimal online EV charging strategy has been studied as in the references above, the online aggregate EV power flexibility characterization problem has not been well explored yet. 
Compared to the online EV charging problem, the online aggregate EV power flexibility characterization is more complicated in three aspects: 1) While the modeling of optimal EV charging problem is clear, how to formulate an optimization model for characterizing aggregate EV power flexibility is an open question. 2) In the online EV charging problem, when the EV charging strategy for the current time slot is determined, its impact on the next time slot (e.g., next time slot's initial EV SoC) can be calculated. However, given the current aggregate EV power flexibility, it is difficult to figure out its future impacts since we do not know which charging strategy within the flexibility region will be chosen. 3) When turning an offline optimal EV charging problem into its online counterpart, we focus on the optimality issue. In contrast, we must verify that the feasible set remains nearly unchanged in the aggregate EV power flexibility characterization.

To fill the research gaps above, this paper proposes a real-time feedback based online aggregate EV power flexibility characterization method.
First, we propose an optimization model to compute the aggregate EV power flexibility region, which is characterized by an upper bound and a lower bound in each time slot. It is theoretically proved that any trajectory within the obtained region is achievable. Second, Lyapunov optimization is leveraged to formulate the online counterpart. Charging task queues and delay-aware virtual queues are constructed to ensure the satisfaction of charging requirements and to handle uncertainties. With the drift-plus-penalty technique, a balance between queues stability and flexibility region value maximization is achieved. Furthermore, we use real-time dispatch information as feedback to improve the accuracy of the aggregate power flexibility region characterization. Particularly, the aggregator updates the virtual queues based on the actual dispatch feedback by the system operator. 
Our main contributions are two-fold:

1) \textbf{Model}. We propose an offline optimization model to characterize the aggregate EV power flexibility region. It decomposes the time-coupled flexibility region into each time slot and gives their lower and upper bounds. Compared to most prior works that provide an outer approximation of the flexibility region, we prove that any trajectory within the region obtained by our model is achievable so that feasibility can be guaranteed. Then, by categorizing the EVs according to their allowed charging delays, we develop a counterpart of the offline model that enables the further utilization of the Lyapunov optimization framework. 

 2) \textbf{Algorithm}. A real-time feedback based online algorithm is developed to derive the aggregate EV power flexibility region. First, to fit into the Lyapunov optimization framework, charging task queues and delay-aware virtual queues are introduced to reformulate the model. Then, a drift-plus-penalty term is constructed and by minimizing its upper bound, an online algorithm is developed. We prove that the charging delays for EVs will not exceed their maximum allowed values even though such constraints are not explicitly considered in the algorithm. The bound of optimality gap between the offline and online outcomes is provided theoretically. Furthermore, real-time dispatch strategy-based feedback is designed and integrated into the online algorithm. Compared to offline \cite{mouli2016system,yang2018distributed,wu2017optimal,wu2018stochastic,yan2018optimized} and online MPC-based \cite{2019zheng} methods, the proposed real-time feedback based online algorithm is prediction-free and can adapt to uncertainties such as random electricity prices and EV charging behaviors. Furthermore, it can make use of the most recent information, allowing it to even outperform the offline model with full knowledge of future uncertainty realizations but without the updated dispatch information.


To summarize, the proposed method is more practical than existing methods in the following ways: 
1) \emph{It runs in a prediction-free manner}, which distinguishes it from the offline and the MPC-based methods that require prior knowledge of future information. Relying solely on the up-to-date information and states, the proposed method can better adapt to uncertainties like random electricity prices and EV charging behaviors. 
2) \emph{It can utilize the most recent dispatch information}, which allows it to even outperform the offline method. This is because that the offline method, though assuming complete future EV and electricity price information, cannot account for the impact of dispatch decisions on future aggregate EV power flexibility. Therefore, the aggregate flexibility may be underestimated.
3) \emph{It has higher computational efficiency}. In contrast to the RL methods, which require a large amount of data and time for training before implementation, the proposed method does not require a training process. This makes it more efficient and practical, especially in situations where there is insufficient data.

The rest of this paper is organized as follows. Section \ref{sec:offline} formulates an offline model for deriving the aggregate EV charging power flexibility region. Sections \ref{sec:online} and \ref{sec:feedback} introduce the Lyapunov optimization method and the real-time feedback design, respectively, to generate the flexibility region in an online manner. Simulation results are presented in Section \ref{sec:result}. Finally, Section \ref{sec:conclu} concludes this paper.

\section{Problem Formulation}\label{sec:offline}
In this section, we first introduce the concept of aggregate EV power flexibility and then formulate an offline optimization problem to approximate it.

\subsection{Aggregate EV Power Flexibility Region}
\begin{figure}[htbp]
  \centering
  \includegraphics[width=0.48\textwidth]{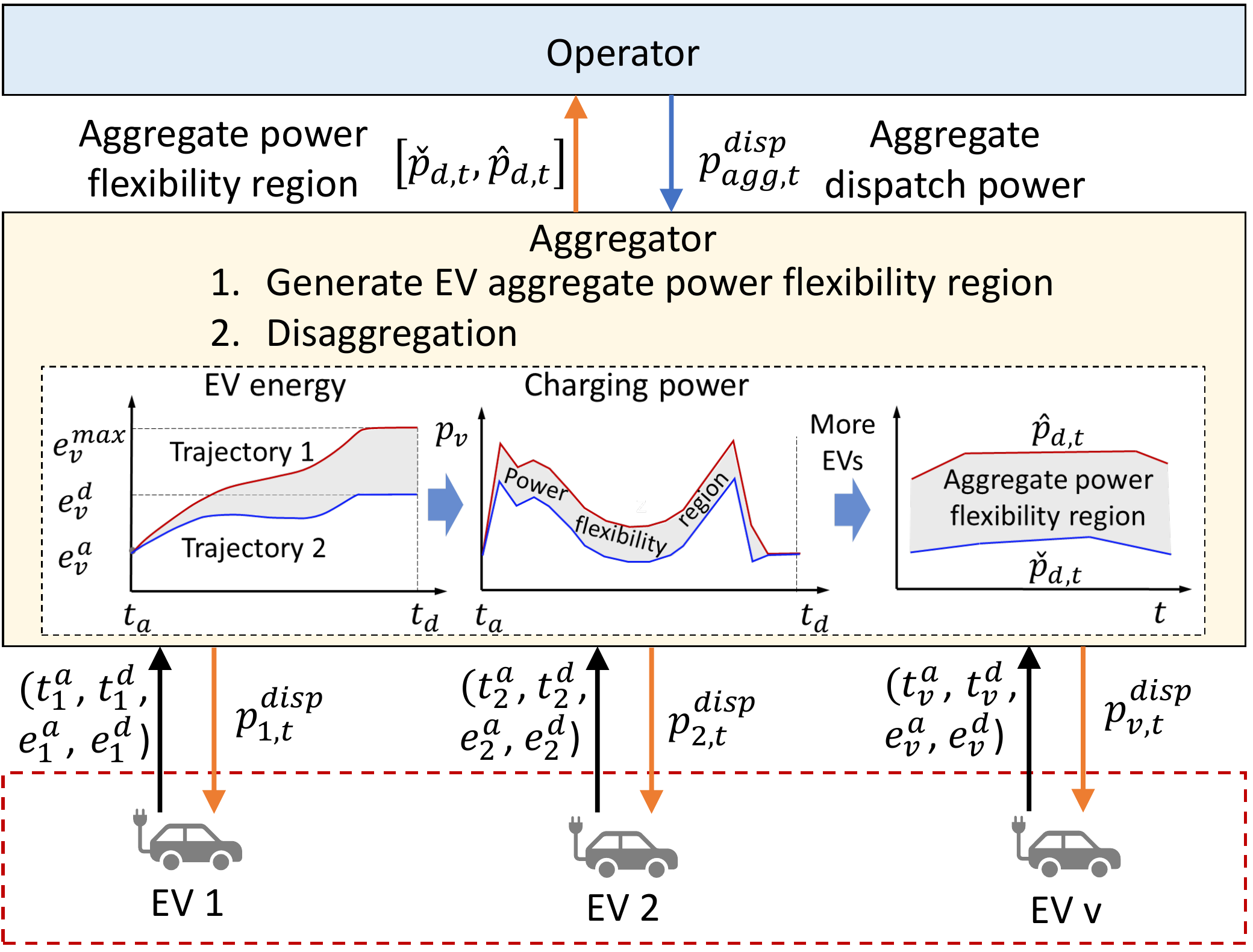}\\
  \caption{System diagram and illustration of EV power flexibility.}\label{fig:sys}
\end{figure}

As shown in Fig. \ref{fig:sys}, when an EV $v\in\mathcal{V}$ arrives at the charging station, it submits its charging task to the aggregator. The task is described by $(t_v^a, t_v^d, e_v^{a}, e_v^{d})$, where $t_v^a$ is its arrival time, $t_v^d$ is its departure time, $e_v^{a}$ is the initial battery energy level at $t_v^a$, and $e_v^{d}$ is the desired energy level when it leaves. For the EV $v$, the maximum allowed charging delay is $R_v=t_v^d-t_v^a$. The EV charging task needs to be completed within this declared time duration.
With the submitted information, the aggregator can flexibly charge the EVs to meet the charging requirements. Two possible trajectories to meet the EV charging needs are depicted in Fig. \ref{fig:sys}. Let $\{p_{v,t}^c,\forall t\}$ be the charging power of EV $v$ over time. The range that the charging power can vary within is called the \emph{power flexibility region} of EV $v$. Similary, the range that $\{\sum_v p_{v,t}^c ,\forall t\}$ can vary within is called the \emph{aggregate EV power flexibility region} of the charging station.

However, it is difficult to characterize the EV power flexibility for each time slot due to the temporal-coupled EV charging constraints. The EV power flexibility in the current time slot is affected by those in the past time slots and further affects those in the future time slots. This is different from the traditional controllable generators whose flexibility can be described by the minimum and maximum power outputs in each time slot. In the following, we aim to derive an aggregate EV power flexibility region that: 1) is time-decoupled so that it can be used in real-time power system operation; and 2) any trajectory within the region can meet the EV charging requirements.

\subsection{Offline Model}
Suppose there are $T$ time slots, indexed by $t \in \mathcal{T}=\{1,...,T\}$. For each time slot $t\in\mathcal{T}$, the desired time-decoupled aggregate EV power flexibility region is given by an interval $[{\check{p}}_{d,t},{\hat{p}}_{d,t}]$ that the total EV charging power can vary within.
The aggregate EV power flexibility region over time is composed of a series of intervals for each time slot, i.e., $[{\check{p}}_{d,1},{\hat{p}}_{d,1}]\times \cdots \times [{\check{p}}_{d,T},{\hat{p}}_{d,T}]$.
These intervals can be specified by a lower power trajectory $\{\check p_{d,t},\forall t\}$ and an upper power trajectory $\{\hat p_{d,t},\forall t\}$. To obtain the lower and upper power trajectories, we formulate the following offline optimization problem: 
\begin{subequations}
\label{eq:flexibility}
\begin{equation}
\textbf{P1:}~\max \limits_{\hat{p}_{d,t},\check{p}_{d,t},\forall t} \mathcal{F}:=\lim\limits_{T\rightarrow\infty}\frac{1}{T}\sum\limits_{t=1}^{T}\mathbb{E}\Big[F_{t}\Big],\label{equ:drobj}
\end{equation}
where
\begin{gather}
F_{t}=\pi_t(\hat{p}_{d,t}-\check{p}_{d,t}),\forall t,\label{eq:Ft}
\end{gather}
subject to
\begin{gather}
  \hat{p}_{d,t}=\sum_{v\in\mathcal{V}}\hat{p}^c_{v,t},\forall t, \label{equ:ubpdi}\\
  0\leq \hat{p}_{v,t}^{c}\leq p_{v}^{max},\forall v, \forall t\in[t_{v}^a,t_{v}^d],\label{equ:ubpd}\\
  \hat{e}_{v,t+1}= \hat{e}_{v,t}+\eta_c\hat{p}_{v,t}^{c}\Delta t,\forall v,\forall t \ne T,\label{equ:ube}\\
  \hat{e}_{v,t_v^a}= e_v^{a},~\hat{e}_{v,t_v^d}\geq e_v^{d},\forall v,\label{equ:ubsntatd}\\
  e^{min}_v\leq \hat{e}_{v,t}\leq e^{max}_v,\forall v,\forall t, \label{equ:uberanges}\\
  \check{p}_{d,t}=\sum_{v\in\mathcal{V}}\check{p}^c_{v,t},\forall t,\label{equ:lbpdi}\\
  0\leq \check{p}_{v,t}^{c}\leq p_{v}^{max},\forall v,\forall t\in[t_{v}^a,t_{v}^d],\\
  \check{e}_{v,t+1}= \check{e}_{v,t}+\eta_c\check{p}_{v,t}^{c}\Delta t,\forall v,\forall t \ne T,\label{equ:lbe}\\
  \check{e}_{v,t_v^a}= e_v^{a},~\check{e}_{v,t_v^d}\geq e_v^{d},\forall v,\\
  e^{min}_v\leq \check{e}_{v,t}\leq e^{max}_v,\forall v,\forall t,\label{equ:lberanges}\\
  \check p_{d,t} \leq \hat p_{d,t},\forall t, \label{equ:joint}\\
  \hat{p}^c_{v,t}=0, \check{p}^c_{v,t}=0, \forall v,\forall t \notin[t_{v}^a,t_{v}^d] \label{eq:evava}
\end{gather}
\end{subequations}

The problem \textbf{P1} aims to maximize the total value of aggregate EV power flexibility. In the objective function (\ref{equ:drobj})-\eqref{eq:Ft}, $\pi_t,\forall t$ are the electricity prices, showing the unit value of EV power flexibility at different time slots. 
Constraint (\ref{equ:ubpdi}) defines the upper bound of the aggregate EV power flexibility region. The charging power of an EV $v$ is limited by (\ref{equ:ubpd}), where $p_{v}^{max}$ is the maximum charging power.
Constraint (\ref{equ:ubsntatd}) defines the EV's initial energy level and the charging requirement.
(\ref{equ:ube}) and (\ref{equ:uberanges}) describe the EV energy dynamics and battery capacity.
Similarly, (\ref{equ:lbpdi})-(\ref{equ:lberanges}) are the constraints related to the lower bound of the aggregate EV power flexibility region. \eqref{equ:joint} is the joint constraint ensuring that $\{\hat p_{d,t},\forall t\}$ and $\{\check p_{d,t},\forall t\}$ provide the upper and lower bounds, respectively.
\eqref{eq:evava} restricts the charging of EVs to their declared parking time.

\begin{proposition} \label{prop-1}
Any aggregate EV charging power trajectory within $[\check p_{d,1}, \hat p_{d,1}] \times 
 \cdots \times [\check p_{d,T},\hat p_{d,T}]$ is achievable.    
\end{proposition}

The proof of Proposition \ref{prop-1} can be found in Appendix \ref{appendix-1}. Despite this nice property, the offline problem \textbf{P1} may be impractical since it requires complete knowledge of the future EV charging tasks and future electricity prices. To address this issue, in the next section, we will first propose a more flexible form of \eqref{eq:flexibility} and adopt the Lyapunov optimization framework to turn the offline problem into its online counterpart. 
Furthermore, to take the impact of real-time dispatch decisions on the future aggregate EV power flexibility into account, a real-time feedback based online flexibility characterization method is developed in Section~\ref{sec:feedback}. 

\section{Online Algorithm}\label{sec:online}
In this section, we adopt the Lyapunov optimization framework to solve the offline problem \textbf{P1} in an online manner. The proposed algorithm can output an aggregate EV power flexibility region with an economic value close to \textbf{P1}.

\subsection{Problem Modification}\label{subsec:pm}
As mentioned above, the charging station serves dozens of EVs every day, and each EV arrives along with a charging task, i.e.,  $(t_v^a, t_v^d, e_v^{a}, e_v^{d})$. The EV charging tasks can be first stored in a queue and be served later according to a first-in-first-out basis.
Since different EVs may have different allowed charging delays, we use multiple queues to classify and collect the EV charging tasks. Suppose there are $G$ types of allowed charging delays $R_g$s, each of which is indexed by $g\in\mathcal{G}=\{1,2,...,G\}$. We construct $G$ queues to collect the respective charging tasks, and each queue is denoted by $Q_g$.
For queue $Q_g$, $Q_{g,t}$ refers to its charging task backlog in time slot $t$. The queue backlog growth is described by
\begin{equation}\label{eq:Qgt}
    Q_{g,t+1}=\max[Q_{g,t}-x_{g,t},0]+a_{g,t},
\end{equation}
where $x_{g,t}$ is the charging power for EVs in group $g$ at time $t$, and $a_{g,t}$ is the arrival rate of EV charging tasks of group $g$ at time $t$. In particular, $a_{g,t}$ is the sum of the energy demand of all EVs in group $g$ that arrive at the beginning of time $t$, i.e.,
\begin{equation}
    a_{g,t} = \sum_{v \in \mathcal{V}_g}a_{v,t},
\end{equation}
where $a_{v,t}$ is the charging demand of EV $v$ in time slot $t$. $\mathcal{V}_g$ is the set of EVs in group $g$.

Recalling that our target in \textbf{P1} is to derive an upper bound and a lower bound for the aggregate EV power flexibility region, we correspondingly define the upper bound queue $\hat{Q}_{g,t}$ and the lower bound queue $\check{Q}_{g,t}$.
Similar to \eqref{eq:Qgt}, we have
\begin{gather}
    \hat{Q}_{g,t+1}=\max[\hat{Q}_{g,t}-\hat{x}_{g,t},0]+\hat{a}_{g,t},\label{eq:Qgt-ub}\\
    \check{Q}_{g,t+1}=\max[\check{Q}_{g,t}-\check{x}_{g,t},0]+\check{a}_{g,t},\label{eq:Qgt-lb}
\end{gather}
where $\hat{x}_{g,t}$ and $\check{x}_{g,t}$ are the charging power for the upper and lower bound queues, respectively, i.e., $\hat{x}_{g,t}=\sum_{v \in \mathcal{V}_g} \hat p_{v,t}^c$ and $\check{x}_{g,t}=\sum_{v \in \mathcal{V}_g} \check p_{v,t}^c$.

The upper and lower bounds of arriving charging demand, i.e., $\hat{a}_{g,t}$ and $\check{a}_{g,t}$, are determined by
\begin{gather}
    \hat{a}_{g,t} = \sum_{v\in \mathcal{V}_g}\hat{a}_{v,t},~
    \check{a}_{g,t} = \sum_{v \in \mathcal{V}_g}\check{a}_{v,t}.
\end{gather}

Particularly, the lower bound of arriving charging demand $\check{a}_{v,t}$ can be determined in the following charging as soon as possible way,
\begin{gather}\label{eq:agvt-lb}
    \check{a}_{v,t}=\left\{
    \begin{array}{ll}
    p_{v}^{max},     & t_{v}^a\leq t<\lfloor \check{t}_{v}^{min}\rfloor+t_{v}^a \\
    \check{e}_{v}^{cha}/\eta_{c}-\lfloor \check{t}_{v}^{min}\rfloor p_{v}^{max},  & t=\lfloor \check{t}_{v}^{min}\rfloor+t_{v}^a \\
    0,           & \text{otherwise}
    \end{array}\right.,
\end{gather}
where $\check{e}_{v}^{cha}=e_{v}^{d}-e_{v}^{a}$, $\check{t}_{v}^{min}$ is the minimum required charging time for EV $v$ determined by $\check{t}_{v}^{min}=\frac{\check{e}_{v}^{cha}}{p_{v}^{max}\eta_{c}}$, and $\lfloor.\rfloor$ means rounding down to the nearest integer.

Different from $\check{a}_{v,t}$, the upper bound of arrival charging demand $\hat{a}_{v,t}$ is determined using the maximum charging demand $e_{v}^{max}$ instead of $e_{v}^{d}$.
Let $\hat{e}_{v}^{cha}=e_{v}^{max}-e_{v}^{a}$, $\hat{t}_{v}^{min}=\frac{\hat{e}_{v}^{cha}}{p_{v}^{max}\eta_{c}}$, and then $\hat{a}_{v,t}$ can be determined by
\begin{gather}\label{eq:agvt-ub}
    \hat{a}_{v,t}=\left\{
    \begin{array}{ll}
    p_{v}^{max},     & t_{v}^a\leq t<\lfloor \hat{t}_{v}^{min}\rfloor+t_{v}^a \\
    \hat{e}_{v}^{cha}/\eta_{c}-\lfloor \hat{t}_{v}^{min}\rfloor p_{v}^{max},     & t=\lfloor \hat{t}_{v}^{min}\rfloor+t_{v}^a \\
    0,           & \text{otherwise.}
    \end{array}\right.
\end{gather}

We then formulate the aggregate EV power flexibility characterization problem as follows:

\begin{subequations}\label{eq:flexibility-p2}
\begin{equation}
\textbf{P2:}~\min \limits_{\hat{x}_{g,t},\check{x}_{g,t},\forall t}\lim\limits_{T\rightarrow\infty}\frac{1}{T}\sum\limits_{t=1}^{T}\mathbb{E}\Big[-F_{t}\Big],\label{equ:drobj-p2}
\end{equation}
subject to
\begin{gather}
  \lim\limits_{T\rightarrow\infty} \frac{1}{T}{\sum_{t=1}^T \mathbb{E}[\hat{a}_{g,t}-\hat x_{g,t}]} \le 0 ,\forall g \label{eq:Qgub-lim} \\
    \lim\limits_{T\rightarrow\infty}\frac{1}{T}{\sum_{t=1}^T \mathbb{E}[\check{a}_{g,t}-\check{x}_{g,t}]} \le 0,\forall g \label{eq:Qglb-lim}\\
  0\leq\hat{x}_{g,t}\leq x_{g,max},\forall g,\forall t, \label{eq:xgub}\\
  0\leq\check{x}_{g,t}\leq x_{g,max},\forall g,\forall t, \label{eq:xglb}\\
  \hat{x}_{g,t} \ge \check{x}_{g,t},\forall g,\forall t,  \label{eq:xgublb}
\end{gather}
\end{subequations}
where $x_{g,max}=\sum_{v \in \mathcal{V}_g} p_v^{max}$. Constraint \eqref{eq:Qgub-lim} ensures that if using the upper bound trajectory $\{\hat{x}_{g,t},\forall t\}$, the charging requirement can be satisfied in a time-average sense. Constraint \eqref{eq:Qglb-lim} poses a similar requirement for the lower bound trajectory $\{\check{x}_{g,t},\forall t\}$. 
Constraints \eqref{eq:xgub} and \eqref{eq:xglb} bound the upper and lower trajectories of the aggregate EV power flexibility region of group $g$. The upper bound is no less than the lower bound, as shown in \eqref{eq:xgublb}. The \textbf{P2} provides a counterpart problem for \textbf{P1}. Similar to Proposition \ref{prop-1}, we can prove that any trajectory between $[\check x_{g,1},\hat x_{g,1}] \times ... \times [\check x_{g,T},\hat x_{g,T}]$ is achievable. However, the allowed charging delay is still not considered in \textbf{P2}, which may result in unfulfilled EV charging tasks upon departure. This issue will be addressed in Section \ref{secIII-C}.

Before that, we first relax \eqref{eq:Qgub-lim} and \eqref{eq:Qglb-lim} into the mean-rate-stable constraints on queues $\hat Q_{g,t}$ and $\check Q_{g,t}$, respectively. Based on \eqref{eq:Qgub-lim}, \eqref{eq:Qglb-lim}, and the definitions of $\hat Q_{g,t+1}, \check Q_{g,t+1}$ in \eqref{eq:Qgt-ub}-\eqref{eq:Qgt-lb}, we have 
\begin{align}\label{eq-1}
    \hat Q_{g,t+1} -\hat a_{g,t} \ge \hat Q_{g,t} -\hat x_{g,t},\forall t
\end{align}
Summing \eqref{eq-1} up over all $t$ and divide both sides by $T$ yields
\begin{align}
     \frac{\mathbb{E}[\hat{Q}_{g,T+1}]-\mathbb{E}[\hat{Q}_{g,1}]}{T} \ge \frac{\sum_{t=1}^T \mathbb{E}[\hat{a}_{g,t}-\hat x_{g,t}]}{T}.
\end{align}
When queues $\hat{Q}_{g,t}$ and $\check{Q}_{g,t}$ are mean rate stable, i.e., $\lim\limits_{T\rightarrow\infty}{\mathbb{E}[\hat{Q}_{g,T+1}]}/{T}=0$ and $\lim\limits_{T\rightarrow\infty}{\mathbb{E}[\check{Q}_{g,T+1}]}/{T}=0$, constraints \eqref{eq:Qgub-lim} and \eqref{eq:Qglb-lim} are satisfied. \eqref{eq:Qgub-lim} and \eqref{eq:Qglb-lim} can be relaxed to the mean-rate-stable constraints on queues $\hat Q_{g,t}$ and $\check Q_{g,t}$.

\subsection{Construct Virtual Queues}
\label{secIII-C}
To overcome the aforementioned charging delay issue, we introduce delay-aware virtual queues,
\begin{gather}
    \hat{Z}_{g,t+1} = \max\{\hat{Z}_{g,t}+\frac{\eta_g}{R_g}\mathbb{I}_{\hat{Q}_{g,t}>0}-\hat{x}_{g,t},0\}, \forall g, \forall t \label{eq:zgt-ub}\\
    \check{Z}_{g,t+1} = \max\{\check{Z}_{g,t}+\frac{\eta_g}{R_g}\mathbb{I}_{\check{Q}_{g,t}>0}-\check{x}_{g,t},0\}, \forall g, \forall t \label{eq:zgt-lb}
\end{gather}
where $\mathbb{I}_{\hat{Q}_{g,t}>0}$ and $\mathbb{I}_{\check{Q}_{g,t}>0}$ are indicator functions of $\hat{Q}_{g,t}$ and $\check{Q}_{g,t}$, respectively. They are equal to 1 if there exist unserved charging tasks in the queues, i.e., $\hat{Q}_{g,t}>0$ and $\check{Q}_{g,t}>0$, respectively. Using $\frac{\eta_g}{R_g}$ to times it, the whole term constitutes a penalty to the virtual queue backlog. $\eta_g$ is a user-defined parameter that can influence the growth rate of the virtual queues. For instance, increasing the value of $\eta_g$ leads to a fast queue growth and a larger backlog value, calling for more attention to accelerate the charging process. We prove that, when the virtual queues have finite upper bounds, with a proper $\eta_g$, the charging delay for EVs in group $g$ is bounded.


\begin{proposition}\label{prop-2}
Suppose $\hat{Q}_{g,t}$, $\check{Q}_{g,t}$, $\hat{Z}_{g,t}$, and $\check{Z}_{g,t}$ have finite upper bounds, i.e., $\hat{Q}_{g,t}\leq \hat{Q}_{g,max}$, $\check{Q}_{g,t}\leq \check{Q}_{g,max}$ $\hat{Z}_{g,t}\leq \hat{Z}_{g,max}$, and $\check{Z}_{g,t}\leq \check{Z}_{g,max}$. The charging delay of all EVs in group $g$ is upper bounded by $\hat{\delta}_{g,max}$ and $\check{\delta}_{g,max}$ time slots, where
\begin{gather}
    \hat{\delta}_{g,max}:=\frac{(\hat{Q}_{g,max}+\hat{Z}_{g,max})R_{g}}{\eta_g},\label{eq:dgmax-ub}\\
    \check{\delta}_{g,max}:=\frac{(\check{Q}_{g,max}+\check{Z}_{g,max})R_{g}}{\eta_g}.\label{eq:dgmax-lb}
\end{gather}
\end{proposition}

The proof of Proposition \ref{prop-2} can be found in Appendix \ref{appendix-2}. It ensures that the charging tasks can always be fulfilled within the available charging periods by properly setting the parameters $\eta_g,\forall g$.

\subsection{Lyapunov Optimization}
Based on the charging task queues and delay-aware virtual queues, the Lyapunov optimization framework is applied as follows.

\subsubsection{Lyapunov Function}
First, we define $\boldsymbol{\Theta}_t=(\boldsymbol{\hat{Q}}_{t},\boldsymbol{\hat{Z}}_{t},\boldsymbol{\check{Q}}_{t}, \boldsymbol{\check{Z}}_{t})$ as the concatenated vector of queues, where
\begin{subequations}
\begin{align}
&\boldsymbol{\hat{Q}}_t=(\hat{Q}_{1,t},...,\hat{Q}_{G,t}),\label{eq:qhat-vector}\\
&\boldsymbol{\hat{Z}}_t=(\hat{Z}_{1,t},...,\hat{Z}_{G,t}),\label{eq:zhat-vector}\\
&\boldsymbol{\check{Q}}_t=(\check{Q}_{1,t},...,\check{Q}_{G,t}),\label{eq:qcheck-vector}\\
&\boldsymbol{\check{Z}}_t=(\check{Z}_{1,t},...,\check{Z}_{G,t}).\label{eq:zcheck-vector}
\end{align}
\end{subequations}
The Lyapunov function is then defined as
\begin{equation}\label{equ:LyaFun}
L(\boldsymbol{\Theta}_t)=\frac{1}{2}\sum\limits_{g\in{\mathcal{G}}}\hat{Q}_{g,t}^2 + \frac{1}{2}\sum\limits_{g\in{\mathcal{G}}}\hat{Z}_{g,t}^2+\frac{1}{2}\sum\limits_{g\in{\mathcal{G}}}\check{Q}_{g,t}^2+\frac{1}{2}\sum\limits_{g\in{\mathcal{G}}}\check{Z}_{g,t}^2,
\end{equation}
where $L(\boldsymbol{\Theta}_t)$ can be considered as a measure of the queue size. A smaller $L(\boldsymbol{\Theta}_t)$ is preferred in order to make the (virtual) queues $\hat{Q}_{g,t}$, $\hat{Z}_{g,t}$, $\check{Q}_{g,t}$, and $\check{Z}_{g,t}$ less congested.

\subsubsection{Lyapunov Drift}
The conditional one-time slot Lyapunov drift is defined as follows:
\begin{equation}\label{equ:LyaDrift}
\Delta(\boldsymbol{\Theta}_t)=\mathbb{E}[L(\boldsymbol{\Theta}_{t+1}) - L(\boldsymbol{\Theta}_t)|\boldsymbol{\Theta}_t],
\end{equation}
where the expectation is taken with respect to the random $\boldsymbol{\Theta}_t$.

The Lyapunov drift is a measure of the expectation of the queue size growth given the current state $\boldsymbol{\Theta}_t$.
Intuitively, by minimizing the Lyapunov drift, virtual queues are expected to be stabilized. However, only minimizing the Lyapunov drift may lead to a low aggregate EV power flexibility value. Therefore, we integrate the expected aggregate flexibility value (\ref{eq:Ft}) for the time slot $t$ with (\ref{equ:LyaDrift}).
The drift-plus-penalty term is obtained, i.e.,
\begin{equation}\label{equ:driftPlusP}
\Delta(\boldsymbol{\Theta}_t)+V\mathbb{E}[-F_t|\boldsymbol{\Theta}_t],
\end{equation}
where $V$ is a weight parameter that controls the trade-off between (virtual) queues stability and aggregate EV power flexibility maximization.

\subsubsection{Minimizing the Upper Bound}
\eqref{equ:driftPlusP} is still time-coupled due to the definition of $\Delta(\boldsymbol{\Theta}_t)$. To adapt to online implementation, instead of directly minimizing the drift-plus-penalty term, we minimize the upper bound of \eqref{equ:driftPlusP} to obtain the upper and lower bounds of the aggregate EV power flexibility region. We first calculate the one-time slot Lyapunov drift:
\begin{align}
&L(\boldsymbol{\Theta}_{t+1}) -L(\boldsymbol{\Theta}_t)\nonumber\\ &=\frac{1}{2}\sum\limits_{g\in{\mathcal{G}}}\Big\{\left[\hat{Q}_{g,t+1}^2-\hat{Q}_{g,t}^2\right]+\left[\hat{Z}_{g,t+1}^2-\hat{Z}_{g,t}^2\right]\nonumber \\
&~~~+\left[\check{Q}_{g,t+1}^2-\check{Q}_{g,t}^2\right]+\left[\check{Z}_{g,t+1}^2-\check{Z}_{g,t}^2\right]\Big\}.\label{equ:lyaD}
\end{align}
Use queue $\hat{Q}_{g,t}$ as an example, based on the queue update equation \eqref{eq:Qgt-ub}, we have
\begin{align}
\hat{Q}_{g,t+1}^2 &= \{\max[\hat Q_{g,t}-\hat x_{g,t},0]+\hat a_{g,t}\}^2\nonumber\\
                  &\leq \hat Q_{g,t}^2+\hat{a}_{g,max}^2+\hat{x}_{g,max}^2+2\hat Q_{g,t}(\hat a_{g,t}-\hat x_{g,t}).
\end{align}
Thus,
\begin{align}
\frac{1}{2}\left[\hat{Q}_{g,t+1}^2-\hat{Q}_{g,t}^2\right]\leq \frac{1}{2}\left(\hat{x}_{g,max}^2+\hat{a}_{g,max}^2\right) \nonumber\\
+\hat{Q}_{g,t}\left(\hat{a}_{g,t}-\hat{x}_{g,t}\right).\label{eq:Q2-Ub}
\end{align}
Similarly, for queues $\check{Q}_{g,t}$, $\hat{Z}_{g,t}$, and $\check{Z}_{g,t}$, we have
\begin{align}
\frac{1}{2}\left[\check{Q}_{g,t+1}^2-\check{Q}_{g,t}^2\right]\leq \frac{1}{2}\left(\check{x}_{g,max}^2+\check{a}_{g,max}^2\right)\nonumber\\
+\check{Q}_{g,t}\left(\check{a}_{g,t}-\check{x}_{g,t}\right).\label{eq:Q2-lb}
\end{align}
\begin{align}
\frac{1}{2}[\hat{Z}_{g,t+1}^2-\hat{Z}_{g,t}^2]\leq \frac{1}{2}\max[(\frac{\eta_g}{R_g})^2,\hat{x}_{g,max}^2]\nonumber\\
+\hat{Z}_{g,t}[\frac{\eta_g}{R_g}-\hat{x}_{g,t}].\label{eq:z2-ub}
\end{align}
\begin{align}
\frac{1}{2}[\check{Z}_{g,t+1}^2-\check{Z}_{g,t}^2]\leq \frac{1}{2}\max[(\frac{\eta_g}{R_g})^2,\check{x}_{g,max}^2]\nonumber\\
+\check{Z}_{g,t}[\frac{\eta_g}{R_g}-\check{x}_{g,t}].\label{eq:z2-lb}
\end{align}

We then substitute the inequalities \eqref{eq:Q2-Ub},\eqref{eq:Q2-lb}, \eqref{eq:z2-ub} and \eqref{eq:z2-lb} into the drift-plus-penalty term and yield
\begin{equation}\label{equ:dppInequ}
    \begin{aligned}
&\Delta(\boldsymbol{\Theta}_t)+V\mathbb{E}[-F_t|\boldsymbol{\Theta}_t]\\
&\leq A+V\mathbb{E}[-F_t|\boldsymbol{\Theta}_t]+\sum\limits_{g\in{\mathcal{G}}}\hat{Q}_{g,t}\mathbb{E}\left[\hat{a}_{g,t}-\hat{x}_{g,t}|\boldsymbol{\Theta}_t\right]\\
&+\sum\limits_{g\in{\mathcal{G}}}\check{Q}_{g,t}\mathbb{E}\left[\check{a}_{g,t}-\check{x}_{g,t}|\boldsymbol{\Theta}_t\right]+\sum\limits_{g\in{\mathcal{G}}}\hat{Z}_{g,t}\mathbb{E}\left[-\hat{x}_{g,t}|\boldsymbol{\Theta}_t\right]\\
&+\sum\limits_{g\in{\mathcal{G}}}\check{Z}_{g,t}\mathbb{E}\left[-\check{x}_{g,t}|\boldsymbol{\Theta}_t\right],
\end{aligned}
\end{equation}
where $A$ is a constant, i.e.,
\begin{align*}
A=~ & \frac{1}{2}\sum\limits_{g\in{\mathcal{G}}}(\hat{x}_{g,max}^2+\hat{a}_{g,max}^2)+\frac{1}{2}\sum\limits_{g\in{\mathcal{G}}}\max[(\frac{\eta_g}{R_g})^2,\hat{x}_{g,max}^2]\nonumber\\
~ & +\frac{1}{2}\sum\limits_{g\in{\mathcal{G}}}(\check{x}_{g,max}^2+\check{a}_{g,max}^2)+\frac{1}{2}\sum\limits_{g\in{\mathcal{G}}}\max[(\frac{\eta_g}{R_g})^2,\check{x}_{g,max}^2] \nonumber\\
~ & + \sum_{g \in \mathcal{G}} [\hat Z_{g,max} \frac{\eta_g}{R_g}] + \sum_{g \in \mathcal{G}} [\check Z_{g,max} \frac{\eta_g}{R_g}].
\end{align*}

By reorganizing the expression in \eqref{equ:dppInequ} and ignoring the constant terms, we can obtain the following online optimization problem:
\begin{align}
  \mathbf{P3}: ~\min\limits_{\hat{x}_{g,t},\check{x}_{g,t},\forall g}~ &\sum\limits_{g\in{\mathcal{G}}}(-V\pi_t-\hat{Q}_{g,t}-\hat{Z}_{g,t})\hat{x}_{g,t}\nonumber\\
&+\sum\limits_{g\in{\mathcal{G}}}(V\pi_t-\check{Q}_{g,t}-\check{Z}_{g,t})\check{x}_{g,t},\label{equ:P3}\\
  \hbox{s.t.}~ &\eqref{eq:xgub}-\eqref{eq:xgublb},\nonumber
\end{align}
where $\hat{Q}_{g,t}$, $\check{Q}_{g,t}$, $\check{Z}_{g,t}$, and $\check{Z}_{g,t}$ are first updated based on \eqref{eq:Qgt-ub}, \eqref{eq:Qgt-lb}, \eqref{eq:zgt-ub}, and \eqref{eq:zgt-lb} before solving \textbf{P3} in each time slot.

In each time slot $t$, given the current system queue state $\boldsymbol{\Theta}_t$, the proposed method determines the current upper and lower bounds $\hat{x}_{g,t}$ and $\check{x}_{g,t}$ of the aggregate EV power flexibility region by solving problem $\textbf{P3}$. Hence, the original offline optimization problem $\textbf{P2}$ has been decoupled into simple online (real-time) problems.
Moreover, different from the conventional offline method (see model \eqref{eq:flexibility} and model \eqref{eq:flexibility-p2}) that requires complete information about the arrival time of all EVs at the beginning ($t=1$), the proposed online method is prediction-free and the arrival time is known only when EVs actually arrive. In particular, at time $t$ when the EVs arrive, we evaluate the lower and upper bounds of the arriving charging demand based on \eqref{eq:agvt-lb} and \eqref{eq:agvt-ub}, and push the charging demand into the virtual queues ${\check{Q}}_{g,t}, {\hat{Q}}_{g,t}$. With the virtual queues and the current electricity price $\pi_t$, the real-time flexibility characterization model \eqref{equ:P3} can be formulated and solved to determine ${\hat{p}}_{d,t}$ and ${\check{p}}_{d,t}$. Then, the real-time gains $F_t$ by providing flexibility can be calculated.

Since the modified problem \textbf{P3} is slightly different from the offline one \textbf{P2}, an important issue we care about is: what's the gap between the aggregate EV power flexibility value of the online problem \textbf{P3} and the offline problem \textbf{P2}? 
\begin{proposition}\label{prop-3}
Denote the obtained long-term time-average aggregate EV power flexibility value, defined in \eqref{equ:drobj}, of $\textbf{P2}$ and $\textbf{P3}$ by $\mathcal{F}^{off}$ and $\mathcal{F}^{*}$, respectively. We have
\begin{equation}\label{equ:gap}
0 \le -\mathcal{F}^{*} + \mathcal{F}^{off} \leq \frac{1}{V}A,
\end{equation}
where $A$ is a constant defined in (\ref{equ:dppInequ}). 
\end{proposition}

The proof of Proposition \ref{prop-3} can be found in Appendix \ref{appendix-3}. The optimality gap can be controlled by the parameter $V$. A bigger $V$ leads to a smaller optimality gap but increased queue sizes. In contrast, a smaller $V$ value makes the queues more stable but results in a larger optimality gap.

\section{Disaggregation \\and Real-time Feedback Design}\label{sec:feedback}
In each time slot $t$, given the aggregate EV power flexibility region $[\sum_g \check x_{g,t}^*, \sum_g \hat x_{g,t}^*]$, the power system operator can determine the optimal aggregate EV dispatch strategy. This aggregate dispatch strategy needs to be further disaggregated to obtain the control strategy for each EV, which is studied in this section. Moreover, considering that the current dispatch strategy will influence the aggregate EV power flexibility in future time slots, real-time feedback is designed and integrated with the proposed online flexibility characterization method in Section \ref{sec:online}.


\subsection{Disaggregation}
Suppose the dispatch strategy in time slot $t$ is $p_{agg,t}^{disp} \in [\sum_{g}\check{x}_{g,t}^*,\sum_{g}\hat{x}_{g,t}^*]$. This can be determined by the system operator through solving an economic dispatch problem based on the up-to-date information (e.g., the grid-side renewable generation). Since this paper focuses on the online characterization of aggregate EV power flexibility, the economic dispatch problem by the operator is omitted for simplicity. Interested readers can refer to \cite{chen2022robust}. 


Let the dispatch ratio $\alpha_t$ be
\begin{equation}\label{eq:ratio}
    \alpha_t=\frac{{p}_{agg,t}^{disp}-\sum_g \check{x}_{g,t}^*}{\sum_{g}\hat{x}_{g,t}^*-\sum_{g}\check{x}_{g,t}^*}.
\end{equation}
Then, the dispatched power $p_{g,t}^{disp}$ for each group $g$ can be determined according to the ratio, i.e.,
\begin{equation}\label{eq:pgt}
    p_{g,t}^{disp}=(1-\alpha_t)\check{x}_{g,t}^*+\alpha_t\hat{x}_{g,t}^*, 
\end{equation}
which satisfies
\begin{equation*}
    p_{agg,t}^{disp}=\sum_{g \in \mathcal{G}}p_{g,t}^{disp}.
\end{equation*}

The next step is to allocate $p_{g,t}^{disp}$ to the EVs in the group $g$. All EVs in the group $g$ are sorted according to their arrival time. Then, we follow a first-in-first-service principle to allocate the energy; namely, the EV that comes earlier will be charged with the maximum charging power. For an EV in group $g$, i.e., $v \in \mathcal{V}_g$, we have
\begin{align}
    p_{v,t}^{disp}=&\min\left\{p_{g,t}^{disp},p_{v}^{max},\frac{e_{v}^{max}-e_{v,t}}{\Delta t}\right\}, \forall v \in \mathcal{V}_g.\label{eq:disagg}
\end{align}
The third term on the right side is used to ensure that the EV will not exceed its allowed maximum energy level.

After the charging assignment for an earlier EV $v \in \mathcal{V}_g$ is completed, the following update procedures will be executed.
\begin{equation}\label{eq:pagg-update}
    p_{g,t}^{disp} \leftarrow (p_{g,t}^{disp} - p_{v,t}^{disp}),
\end{equation}
which means deducting $p_{v,t}^{disp}$ ($v \in \mathcal{V}_g$) from the total remaining dispatched power $p_{g,t}^{disp}$. 

Then, $p_{g,t}^{disp}$ is allocated to the next earliest-arrived EVs until the aggregate EV charging power is completely assigned. At this time, the disaggregation is finished.

The dispatched power disaggregation algorithm is presented in Algorithm \ref{algo:disagg}.
\begin{algorithm}[htbp]
\caption{EV Dispatched Power Disaggregation}\label{algo:disagg}
\begin{algorithmic}[1]
\STATE{\textbf{Initialization}: aggregate EV dispatched power $p_{agg,t}^{disp}$.}
\STATE{Calculate the dispatched aggregate charging power $p_{g,t}^{disp}$ for each group $g$ using \eqref{eq:pgt}.}
\FOR{Each group $g\in\mathcal{G}$}
\FOR{Each EV $v$ in group $g$}
\IF{the EV is not available for charging}
\STATE{Let EV $v$'s charging power $p_{v,t}^{disp}=0,\forall v$.}
\ELSE
\STATE{Calculate $p_{v,t}^{disp}$ according to \eqref{eq:disagg}.}
\STATE{Update the remaining aggregate power via \eqref{eq:pagg-update}.}
\IF{the updated $p_{g,t}^{disp}=0$}
\STATE Break and return to Step 3.
\ENDIF
\ENDIF
\ENDFOR
\ENDFOR
\end{algorithmic}
\end{algorithm}

\subsection{State Update to Improve Power Flexibility Region}

\begin{figure}[!htbp]
  \centering
  \includegraphics[width=0.45\textwidth]{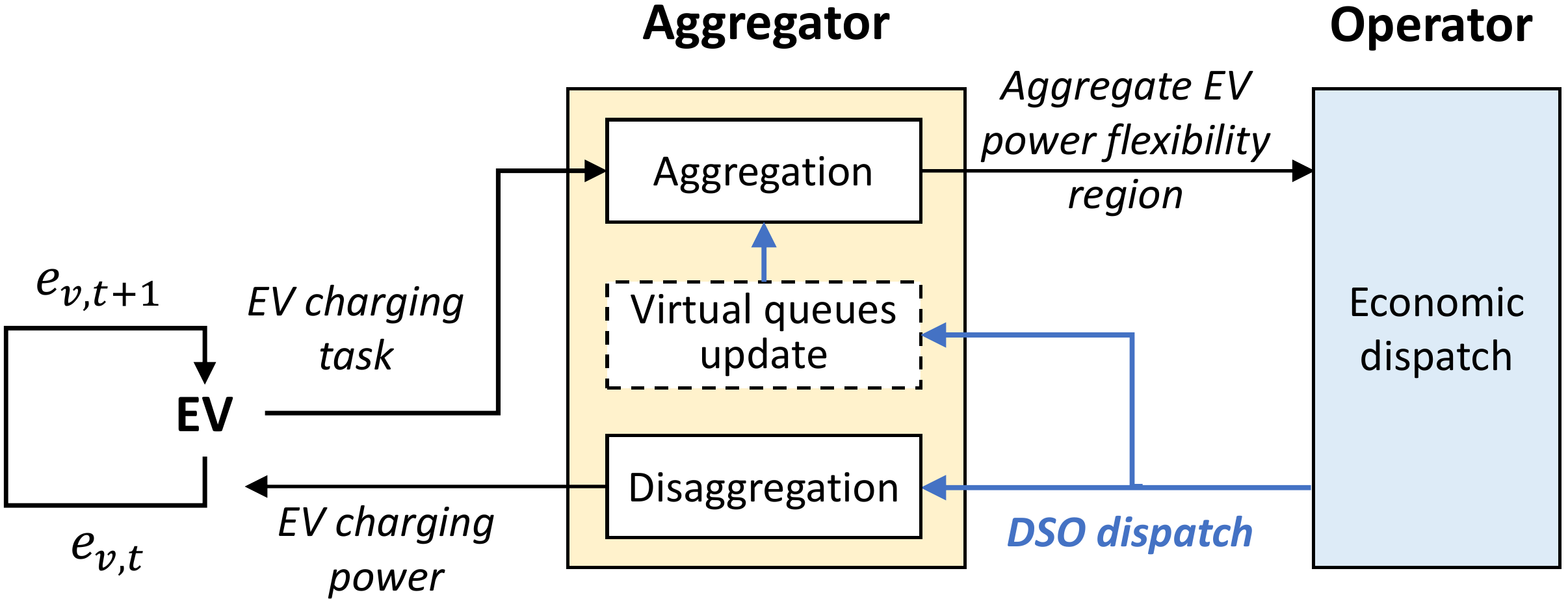}\\
  \caption{Overall procedure of the proposed method with information feedback.}\label{fig:procedure}
\end{figure}

Following the disaggregation procedures in Algorithm 1, we can get the actual EV dispatched charging power $p_{v,t}^{disp},\forall v$. By now, we can move on to the next time slot $t+1$ and evaluate the EV power flexibility by solving problem \textbf{P3}, determine the dispatch strategy, disaggregate the dispatched power, and so on. However, the current actual dispatched EV charging power can affect the future aggregate EV power flexibility, which is ignored in the aforementioned processes.
To address this issue, we develop an online aggregate EV power flexibility characterization method and use the real-time dispatch information as feedback to improve the characterization accuracy.
Note that the feedback in this paper refers to the information feedback for online optimization (like in \cite{li2020real}), rather than the state feedback in control problems.
The timeline is shown in Fig. \ref{fig:procedure}. At time $t$, the aggregator solves the online problem \eqref{equ:P3} to determine the real-time aggregate EV power flexibility region. The region is sent to the system operator who then decides on the real-time aggregate EV dispatch strategy. The dispatch strategy is sent back to the aggregator who conducts disaggregation to obtain the dispatch order for individual EVs. Then, to improve the accuracy of flexibility evaluation, instead of moving directly to the flexibility characterization problem for the next time slot, the aggregator updates the virtual queues based on the actual dispatch information feedback by the system operator first.

To be specific, we change the constraints \eqref{eq:xgub}-\eqref{eq:xglb} for time slot $t$ into
\begin{align}\label{eq:replace}
    \hat x_{g,t}= p_{g,t}^{disp},~  \check x_{g,t}= p_{g,t}^{disp}.
\end{align}
Since $p_{g,t}^{disp} \in [\check x_{g,t}^*, \hat x_{g,t}^*]$, after replacing \eqref{eq:xgub}-\eqref{eq:xglb} with \eqref{eq:replace}, the problem \textbf{P3} is still feasible and denote the optimal solution by $\hat x_{g,t}^{update*}= p_{g,t}^{disp}, \check x_{g,t}^{update*}= p_{g,t}^{disp},\forall g$.  With these updated lower and upper bounds, we update the queues $\hat Q_{g,t+1}$, $\check Q_{g,t+1}$, $\hat Z_{g,t+1}$,$\check Z_{g,t+1}$ according to \eqref{eq:Qgt-ub}, \eqref{eq:Qgt-lb}, \eqref{eq:zgt-ub} and \eqref{eq:zgt-lb}, respectively. Then, we move on to solve problem \textbf{P3} for time slot $t+1$ using the updated $\hat Q_{g,t+1}$, $\check Q_{g,t+1}$, $\hat Z_{g,t+1}$,$\check Z_{g,t+1}$.



A completed description of the proposed methods, including a real-time feedback based online aggregate EV power flexibility characterization method and a EV dispatched charging power disaggregation approach, is shown in Algorithm \ref{algo}.

\begin{algorithm}[htbp]
\caption{Real-time Feedback Based Online Aggregate EV Power Flexibility Characterization and Disaggregation}\label{algo}
\begin{algorithmic}[1]
\renewcommand{\algorithmicrequire}{ \textit{\textbf{Initialization}}}
\REQUIRE
\STATE{Set $t=1$, queues $\hat{Q}_{g,1}=0,\check{Q}_{g,1}=0,\hat{Z}_{g,1}=0,\hat{Z}_{g,1}=0$, parameters $V>0,\eta_g>0$.}
\renewcommand{\algorithmicrequire}{ \textit{\textbf{I. Aggregation}}}
\REQUIRE
\STATE{Aggregator classifies the arriving EVs and pushes them into different queues $\hat{Q}_{g,t}$, $\check{Q}_{g,t}$, $\hat{Z}_{g}$ and $\check{Z}_{g}$ according to their declared allowed charging delays $R_g$.}
\STATE{Solve problem \textbf{P3} to obtain the aggregate EV power flexibility region $[\check{x}_{g,t}^*,\hat{x}_{g,t}^*]$.}
\STATE{Update queues $\hat{Q}_{g,t+1}$, $\check{Q}_{g,t+1}$, $\hat{Z}_{g,t+1}$, and $\check{Z}_{g,t+1}$ according to \eqref{eq:Qgt-ub}, \eqref{eq:Qgt-lb}, \eqref{eq:zgt-ub} and \eqref{eq:zgt-lb}, respectively.}
\renewcommand{\algorithmicrequire}{ \textit{\textbf{II. Dispatch and Disaggregation}}}
\REQUIRE
\STATE{Receive the dispatch decision from the operator.}
\STATE{Perform EV dispatched power disaggregation  according to Algorithm \ref{algo:disagg}.}
\renewcommand{\algorithmicrequire}{ \textit{\textbf{III. Real-Time Feedback and Update}}}
\REQUIRE
\STATE{Update the lower and upper bounds $\hat x_{g,t}^{update*}$, $\check x_{g,t}^{update*},\forall g$ of the EV aggregate power flexibility region.}
\STATE{Update queues $\hat{Q}_{g,t+1}$, $\check{Q}_{g,t+1}$, $\hat{Z}_{g,t+1}$, and $\check{Z}_{g,t+1}$ according to \eqref{eq:Qgt-ub}, \eqref{eq:Qgt-lb}, \eqref{eq:zgt-ub} and \eqref{eq:zgt-lb}, respectively.}
\STATE{Move to the next time slot $t \leftarrow t+1$, and repeat the above steps I-III.}
\end{algorithmic}
\end{algorithm}

\subsection{Discussions on Practical Issues}
In the following, We further discuss some practical issues regarding the proposed model and method below:

(1) \emph{Compatibility with existing electricity markets}.
In the real-time electricity market, participants submit real-time supply-offers and demand-bids to the operator, in the form of price-quantity curves for the next time period. Based on the supply-offers and demand-bids, the operator clears the real-time market to determine the energy quantities and price at equilibrium \cite{widergren2022dso+}. In some countries, the electricity market is still in an early stage, where demand-side participants are treated as price-takers and can only submit the electricity quantities without prices \cite{guo2020power}. The proposed model could be suitable for this kind of electricity markets.
Here, we briefly introduce how the proposed model fits into the real-time electricity market: With the up-to-date electricity price, the EV aggregator solves the problem $\textbf{P3}$ to obtain the lower and upper bounds of the flexibility region, i.e., $[\sum_{g}{\check{x}}_{g,t}^\ast,\sum_{g}{\hat{x}}_{g,t}^\ast]$. Then, the aggregator submits it to the system operator. Based on the received aggregator demand-bids, the operator clears the real-time market and computes the cleared/dispatched energy price and quantities to send back to the EV aggregator.

\begin{figure}[!htbp]
  \centering
  \includegraphics[width=0.45\textwidth]{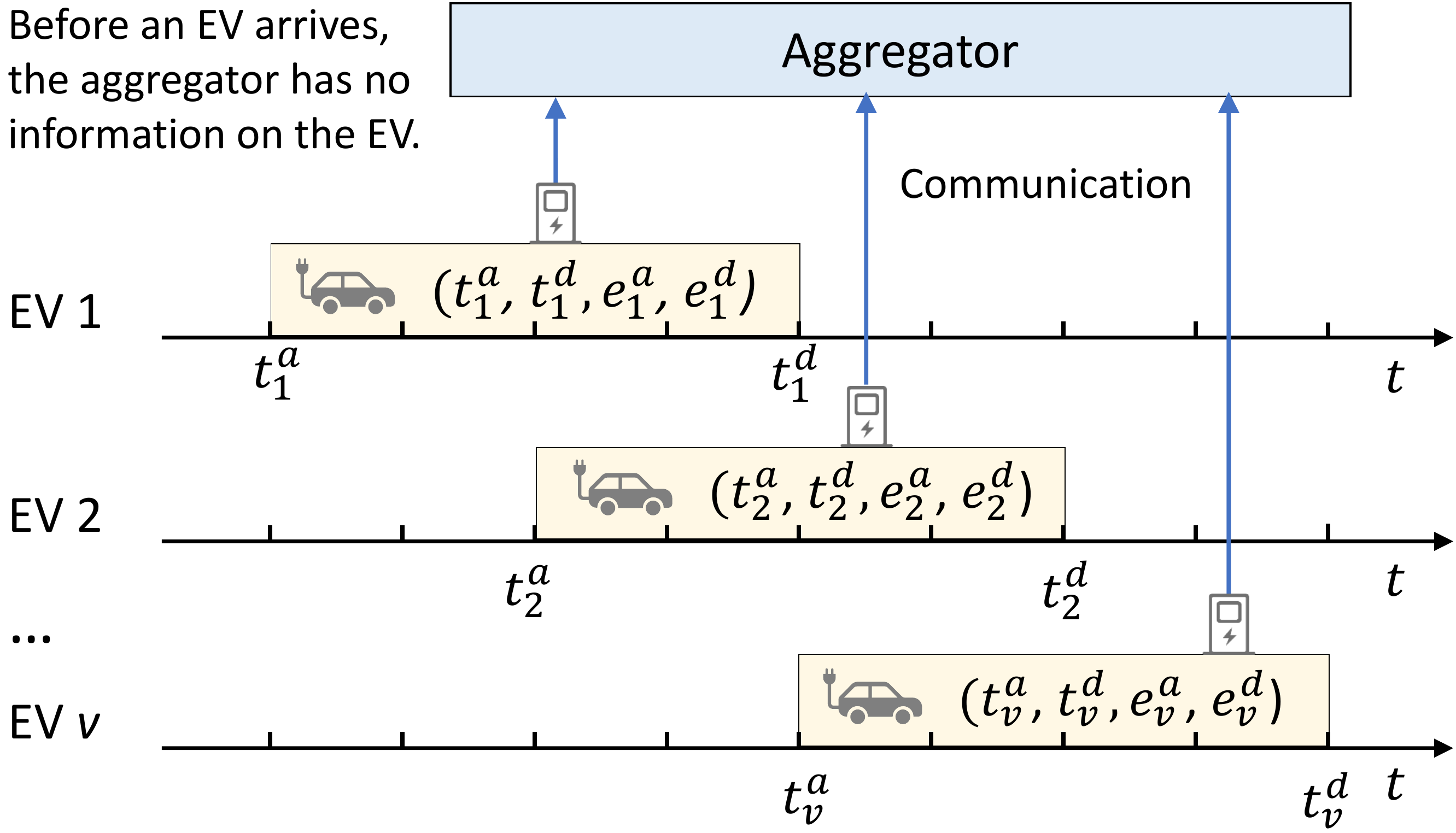}\\
  \caption{Illustration of information exchanged between EVs and the aggregator.}\label{fig:info}
\end{figure}

(2) \emph{Dealing with uncertainties}.
Generally, there are two types of approaches to dealing with uncertainties, i.e., one based on stochastic/robust optimization techniques and the other utilizes online optimization methods \cite{guo2021real}. This paper focuses on the second type of approach by proposing an online flexibility characterization method. The proposed online method is prediction-free and can adapt to different realizations of uncertainties. In particular, the EV related parameters (e.g., the arriving time, the initial battery energy level) are not known in advance, but are revealed only when the EVs actually arrive, as shown in Fig. \ref{fig:info}. In contrast, the conventional offline models rely on predictions of the EVs at the beginning ($t=1$).
It is worth noting that in the test study, the data are randomly generated. As we need to compare the proposed algorithm with other baselines, we fix the EV parameters to keep the settings the same. But this does not mean the proposed method is deterministic.

(3) \emph{Privacy}. In the proposed structure, the aggregator acts as an intermediate between EVs and the system operator. This can avoid EVs directly revealing private information to the operator and thus protect their privacy to some extent. The operator only gets the aggregate EV power flexibility region derived by the aggregator but not the detailed information of individual EVs. We have to admit that an opponent might still be able to reverse-engineer the sensitive data from the aggregate flexibility region information exchanged. How to fully protect EVs' privacy will be studied in our future work.

(4) \emph{Model flexibility}. The proposed model is quite flexible and can readily accommodate the evolving landscape of EV proliferation in future scenarios. For instance, as EV charging technology advances, charging duration could be significantly reduced. To deal with this, the time resolution of the proposed model can be refined to smaller intervals. In addition, to improve the utilization rate of charging infrastructure, the proposed model allow EVs to modify their allowed charging delay during the charging process. To illustrate, suppose an EV initially sets its allowed charging delay $R_g$ to 5 hours. Then, after an hour of charging, the EV may decide to reduce its remaining allowed charging delay to 3 hours.  In this scenario, the EV can report to the aggregator with $R_g=3$ hours and the up-to-date SOC information. The system will treat this EV as a new arrival, assign it to the group according to its newly declared $R_g$, and add corresponding new constraints. 

\section{Simulation Results and Discussions}\label{sec:result}
In this section, we evaluate the performance of the proposed online algorithm and compare it with other approaches. Basically, we test the performance of the proposed algorithm under two situations:
1) Less challenging situations where the EVs have longer allowed charging delay (Sections \ref{secV-A} and \ref{secV-B}). This is commonly found in charging stations located in campus and commercial working areas, where EV owners typically arrive at around 9:00 am and depart at around 6:00 pm.
2) More challenging situations where the EVs have shorter allowed charging delays and higher desired SOC (Section \ref{secV-F}). This may be found in charging stations located in a shopping mall.
The formulated offline and online optimization problems are modeled in MATLAB and solved by GUROBI solver. All simulations are run on a laptop with an Intel Core i7 1.8 GHz CPU and 16 GB RAM.

\subsection{System Setup}
\label{secV-A}
\begin{figure}[!htbp]
  \centering
  \includegraphics[width=0.45\textwidth]{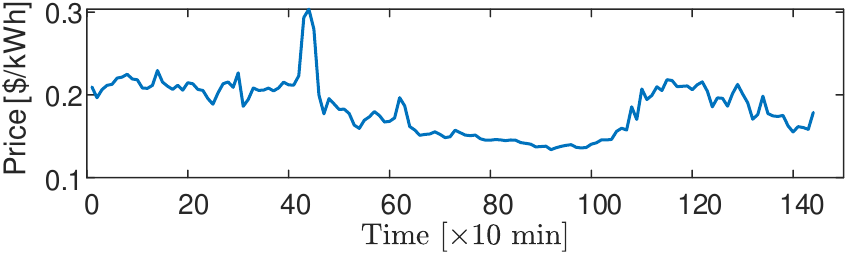}\\
  \caption{Real-time electricity price profile.}\label{fig:price}
\end{figure}

\begin{figure}[!htbp]
  \centering
  \includegraphics[width=0.45\textwidth]{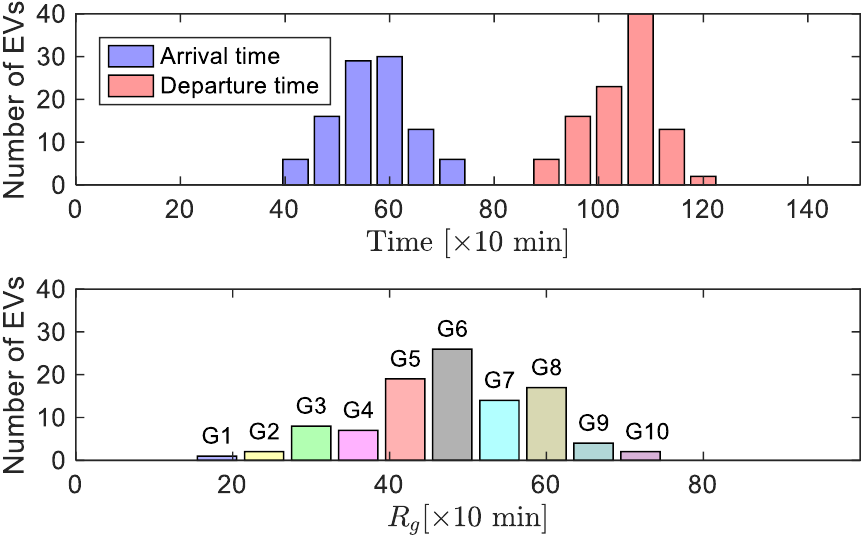}\\
  \caption{Distribution of EV arrival time and departure time, and distribution of $R_g$, where ``G'' means group.}\label{fig:time}
\end{figure}

The entire simulation duration considered is 24 hours, i.e., 144 time slots, each with a time interval of 10 minutes.
To reflect the actual fluctuations in electricity prices, we use the real-time electricity price data obtained from the PJM market \cite{pjm}. The dynamic electricity price data profile is shown in Fig. \ref{fig:price}.
We consider 100 EVs in a commercial area. Their arrival and departure times follow normal distributions $\mathcal{N}(\text{9:00}, (1.2hr)^2)$ and $\mathcal{N}(\text{18:00}, (1.2hr)^2)$ \cite{mohamed2014real}, respectively, as shown in Fig. \ref{fig:time} (top). Based on these, we can calculate the parking time, namely, the allowed charging delay $R_g$, which ranges from 3 to 12 hours. The distribution of $R_g$ is shown in Fig. \ref{fig:time} (bottom). As seen, there are $G=10$ groups, and each group has a different number of EVs.
In addition, the initial battery energy level of each EV is selected from a uniform distribution in $[0.3,0.5] \times$EV battery capacity randomly \cite{jin2014optimized}. We set the required state-of-charge (SOC)\footnote{The SOC of an EV is the ratio between the battery energy level and the battery capacity.} upon departure as 0.5 and the maximum SOC upon departure as 0.9.
To incorporate the diversity of EV charging technical characteristics, we consider 3 commonly-used types of EV battery capacity (24 kWh, 40 kWh, and 60 kWh) and 3 types of charging power (3.3 kW, 6.6 kW, and 10 kW). The weight parameter $V$ is set as 200, and $\eta_g$ is set as $5R_g$ for each group.

\subsection{Effectiveness of the proposed method}
\label{secV-B}
\begin{figure}[!htbp]
  \centering
  \includegraphics[width=0.45\textwidth]{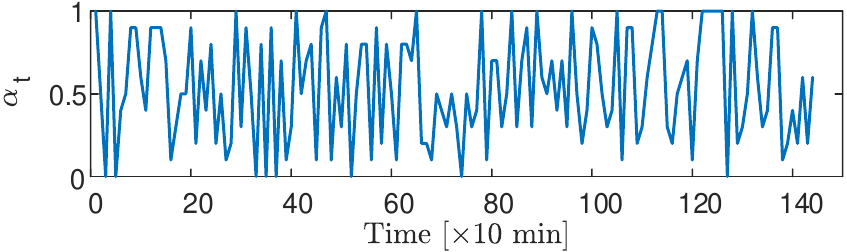}\\
  \caption{Randomly generated dispatch ratio $\alpha_t$.}\label{fig:alpha}
\end{figure}

\begin{figure*}[!htbp]
  \centering
  \includegraphics[width=0.95\textwidth]{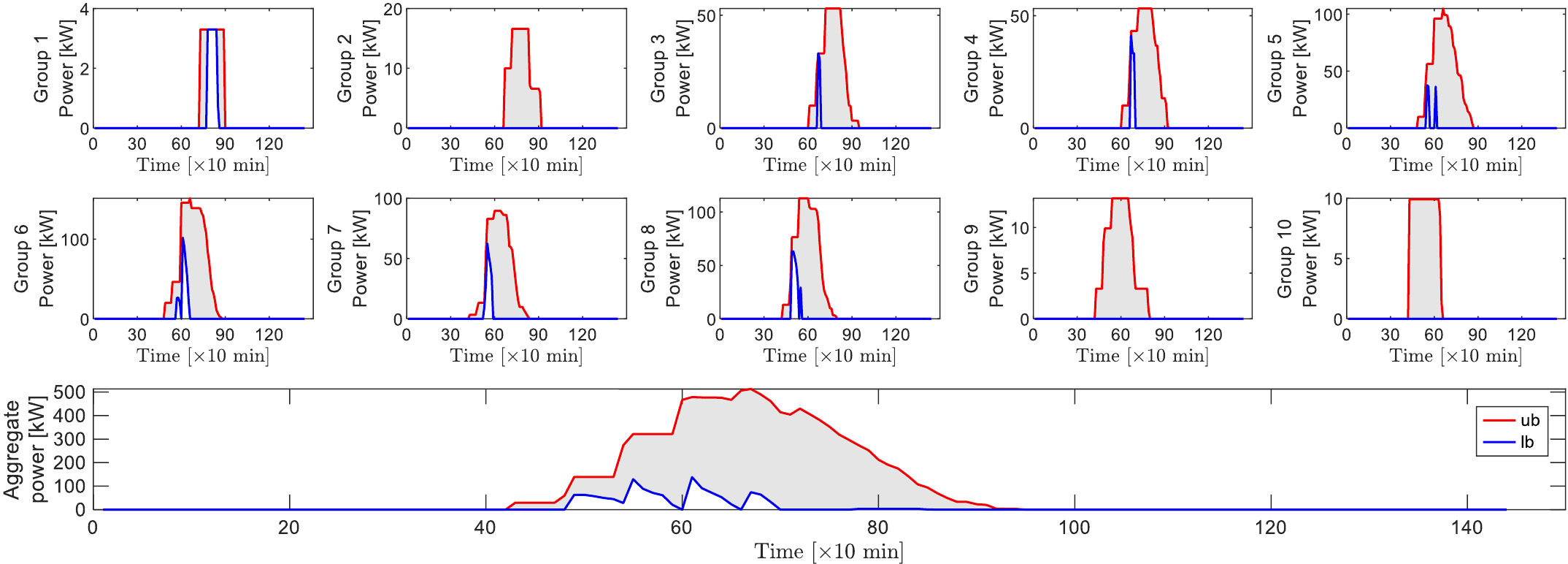}\\
  \caption{The obtained aggregate EV power flexibility region.}\label{fig:flex_prop}
\end{figure*}

We first show what the obtained aggregate EV power flexibility region looks like. Since the power grid dispatch determined by the operator is beyond the scope of this paper, here the dispatch ratio $\alpha_t$ in \eqref{eq:ratio} is randomly generated within the range of $[0,1]$ for each time slot, as shown in Fig. \ref{fig:alpha}.
Based on the generated dispatch ratio, we apply the proposed online flexibility characterization method and real-time feedback in turns (as in Algorithm \ref{algo}) to obtain the aggregate EV power flexibility region (grey area) for each group and the charging station as a whole. The results are shown in Fig. \ref{fig:flex_prop}. As seen, the power flexibility region varies over time. This is because EVs dynamically arrive and leave. 
In addition, we find that the EV power flexibility for groups 1 and 3-8 in some time slots are empty even though the EVs are available for charging. It means that EV charging becomes inelastic in those periods and cannot offer power flexibility. Thanks to the aggregation effect, the aggregate power flexibility region is always non-empty during the available charging time slots.

\begin{figure*}[!htbp]
  \centering
  \includegraphics[width=0.95\textwidth]{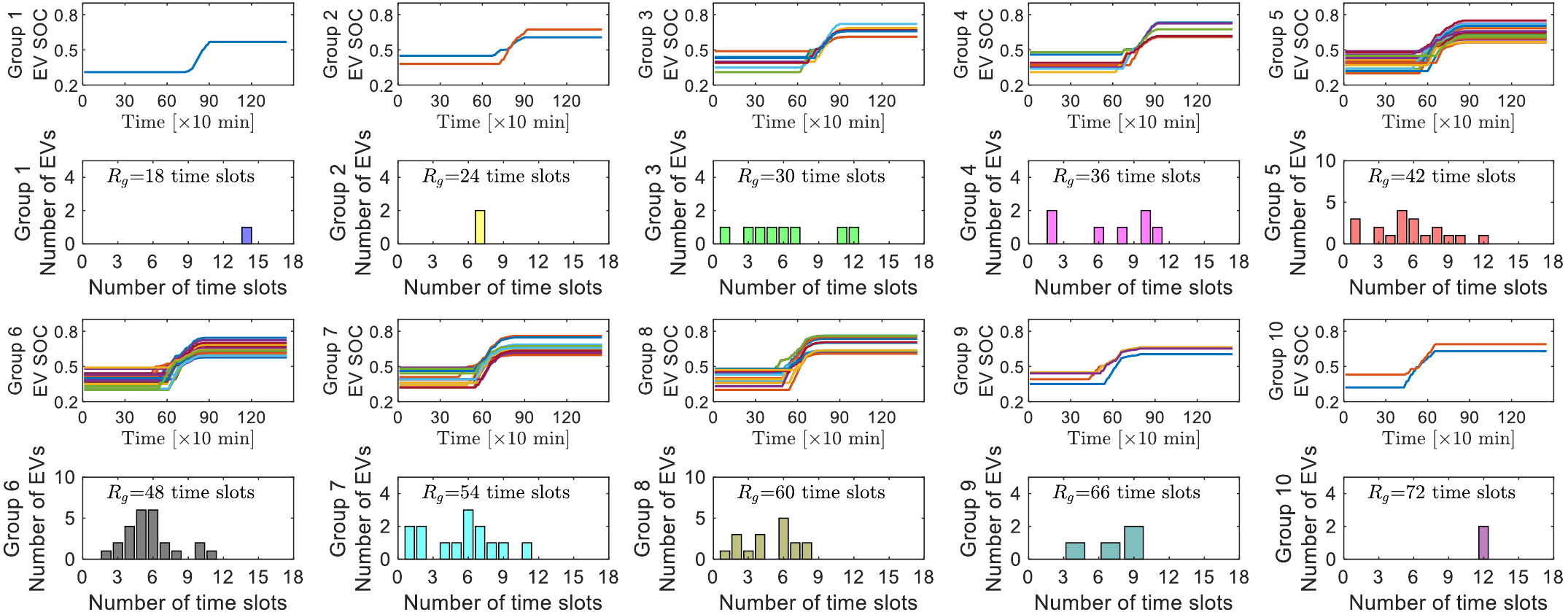}\\
  \caption{EV charging SOC of each group and the number of time slots taken to reach the charging requirement SOC $=0.5$.}\label{fig:soc}
\end{figure*}

To validate the effectiveness of the proposed algorithm, disaggregation of the dispatched EV charging power is performed and we check if the SOC curves of EVs satisfy the charging requirements. Here, if the final EV SOC value can reach or exceed the EV owner's requirement (SOC $\geq 0.5$) upon leaving, then the proposed method is effective. Fig. \ref{fig:soc} shows the actual EV charging SOC curves under the randomly generated dispatch ratio $\alpha_t$ in Fig. \ref{fig:alpha}. Each curve represents an EV. As we can see from the figure, all EVs' final SOC is between 0.5 and 0.7, greater than the required value (0.5) and less than the maximum value (0.9).
Fig. \ref{fig:soc} also shows the number of time slots EVs take to reach the required SoC $=0.5$. We can find that the maximum number of time slots spent is 14 for group 1, 7 for group 2, 12 for group 3, 12 for group 4, 11 for group 5, 11 for group 6, 11 for group 7, 8 for group 8, 9 for group 9, and 12 for group 10. All of them are no more than their respective declared allowed charging delay, i.e., $R_g$. This validates the effectiveness of the proposed algorithm in providing maximum power flexibility while meeting the charging requirements.
 
\begin{figure*}[!htbp]
  \centering
  \includegraphics[width=0.95\textwidth]{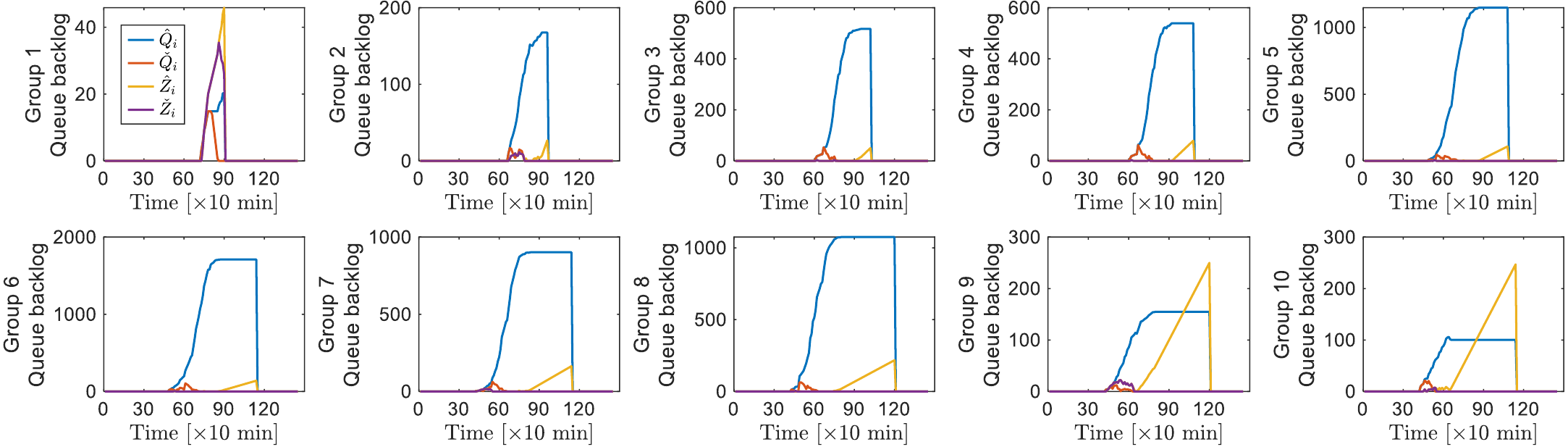}\\
  \caption{Queue backlog of each group.}\label{fig:queue}
\end{figure*}

Furthermore, Fig. \ref{fig:queue} shows the queue backlog evolution of the ten groups over time. Taking group 1 for example, the lower and upper bound queues $\check{Q}_1$ and $\hat{Q}_1$ first increase because EVs arrive and the charging tasks are pushed into the queues. As time moves on, through the EV charging dispatch $p_{v,t}^{disp},\forall v,\forall t$ determined by disaggregation, the EV SOCs gradually increase and reach the minimum charging requirement 0.5. Hence, the lower bound queue $\check Q_1$ becomes zero because the $\check a_{g,t},\forall g,\forall t$ is set using the charging as soon as possible method to meet the minimum charging requirement as in \eqref{eq:agvt-lb}. The upper bound queue $\hat{Q}_1$ is still larger than zero since the EV SOCs have not reached the maximum value 0.9 (see Fig. \ref{fig:soc}), so there is still charging flexibility. 
The changes of delay aware queue $\hat{Z}_1$/$\check{Z}_1$ is influenced by both the charging task queue $\hat{Q}_1$/$\check{Q}_1$ and the actual dispatch charging power $\hat{x}_1$/$\check{x}_1$. Specifically, the initial rapid growth of $\hat{Z}_1$/$\check{Z}_1$ can be attributed to the penalty for nonzero $\hat{Q}_1$/$\check{Q}_1$, which exceeds the dispatch $\hat{x}_1$/$\check{x}_1$. Once $\check{Q}_1$ becomes zero, the dynamics of $\check{Z}_1$ are solely controlled by $\check{x}_1$, and thus $\check{Z}_1$ decreases. 
The queue evolution in other groups can be analyzed similarly.

\emph{Remark}: The proposed algorithm can not only be implemented on MATLAB but also other platforms, such as Python 3.6. The Python program takes 0.7 seconds while the Matlab program takes 2.5 seconds. This demonstrates the computational efficiency of the proposed algorithm regardless of the simulation platform.


\subsection{Performance Evaluation}
\label{secV-C}
To show the advantages of the proposed online algorithm, two widely used benchmarks in the literature are performed.
\begin{itemize}
  \item Benchmark 1 (B1): This is a greedy algorithm that EVs start charging at the maximum charging power upon arrival. Denote the arrival time as $t_0$. When the EV SOC reaches the minimum charging requirement 0.5, the lower bound of charging power $\check p_{v,t}$ turns to zero (time: $t_1$), and the upper bound of charging power $\hat p_{v,t}$ remains to be the maximum charging power until the EV SOC reaches the maximum value 0.9 (time: $t_2$). The aggregate power flexibility region for $[t_0,t_1]$ is empty and for $t\in[t_1,t_2]$ is the region between 0 and the maximum charging power.
  \item Benchmark 2 (B2): This is the offline method. It directly solves \textbf{P1} to obtain the aggregate EV power flexibility regions over the whole time horizon by assuming known future information. Though not realistic, it provides a theoretical benchmark to verify the performance of other methods. But it is worth noting that since it does not take into account the real-time actual dispatch strategy when calculating the aggregate flexibility, its performance may be worse than the proposed real-time feedback based method even though it is an offline method.
\end{itemize}

\begin{figure}[!htbp]
  \centering
  \includegraphics[width=0.35\textwidth]{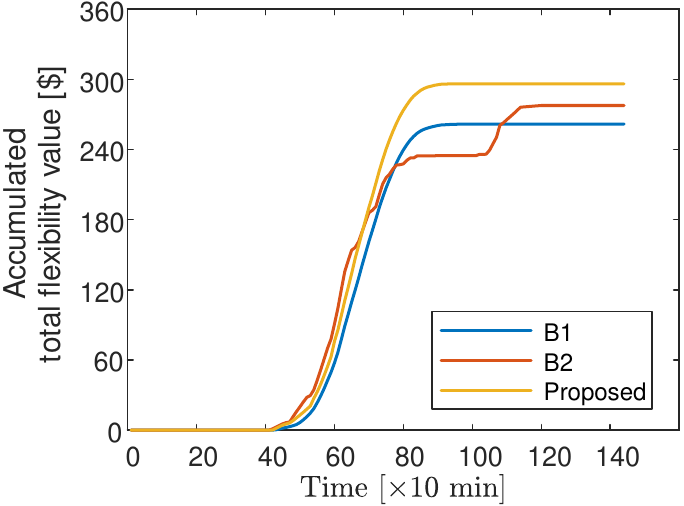}\\
  \caption{Accumulated flexibility value under different methods.}\label{fig:accum}
\end{figure}

Fig. \ref{fig:accum} shows the accumulated flexibility values ($\sum_{\tau=1}^{t}F_{\tau}$) under the three different methods, and TABLE \ref{tab:compare} summarizes the total flexibility value ($\sum_{t=1}^{T}F_t$) under different methods.
The B1, i.e., greedy algorithm, has the worst performance and the lowest total flexibility value due to the myopic strategy.
For B2, since it has complete future knowledge of EV behaviors and real-time electricity prices, it outperforms B1 with an increase of 6.1\%. However, this method is usually impractical since accurate future information is hardly available. Though predictions on future uncertainty realizations may be obtained, the potential prediction errors limit B2's performance.
The proposed online algorithm achieves the best performance with the highest total power flexibility value with an increase of 13.1\%.
This is owing to the fact that the proposed method runs in an online manner with real-time feedback that allows it to utilize the most recent dispatch information to update the state of EVs. In contrast, the offline model is solved at the beginning of the time period ($t=1$) and cannot take into account the real-time dispatch information. In addition, compared to the offline method B2, the proposed method does not require prior knowledge of future information or forecasts, which is more practical.

\begin{table}[!htbp]
  \centering
  \caption{Total flexibility value comparison between B1, B2, and the proposed algorithm (Unit: USD).}\label{tab:compare}
  \begin{tabular}{cccc}
  \hline\hline
    Methods & B1    & B2       & Proposed \\
    \hline
    Value & 261.8  & 277.6      & 296.1 \\
    Improvement & - & 6.1\%   & 13.1\% \\
    \hline\hline
   \end{tabular}
\end{table}

\begin{figure}[!htbp]
  \centering
  \includegraphics[width=0.45\textwidth]{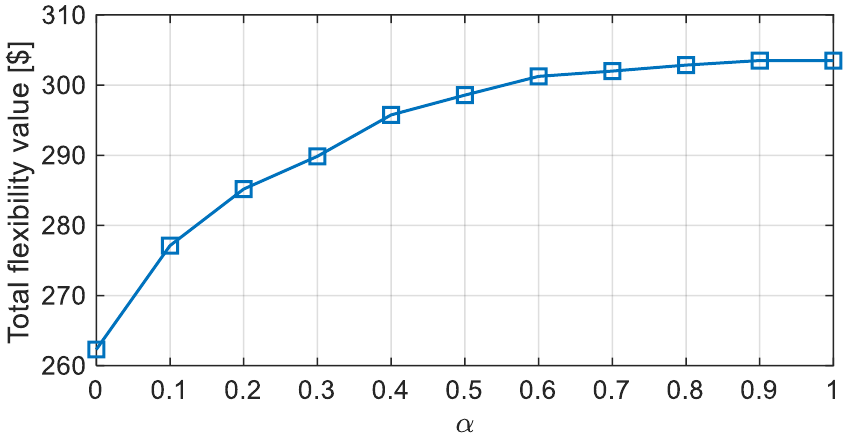} \\
  \caption{The impact of $\alpha$ on total power flexibility value.}
  \label{fig:diff_alpha}
\end{figure}

\begin{figure}[!htbp]
  \centering
  \includegraphics[width=0.45\textwidth]{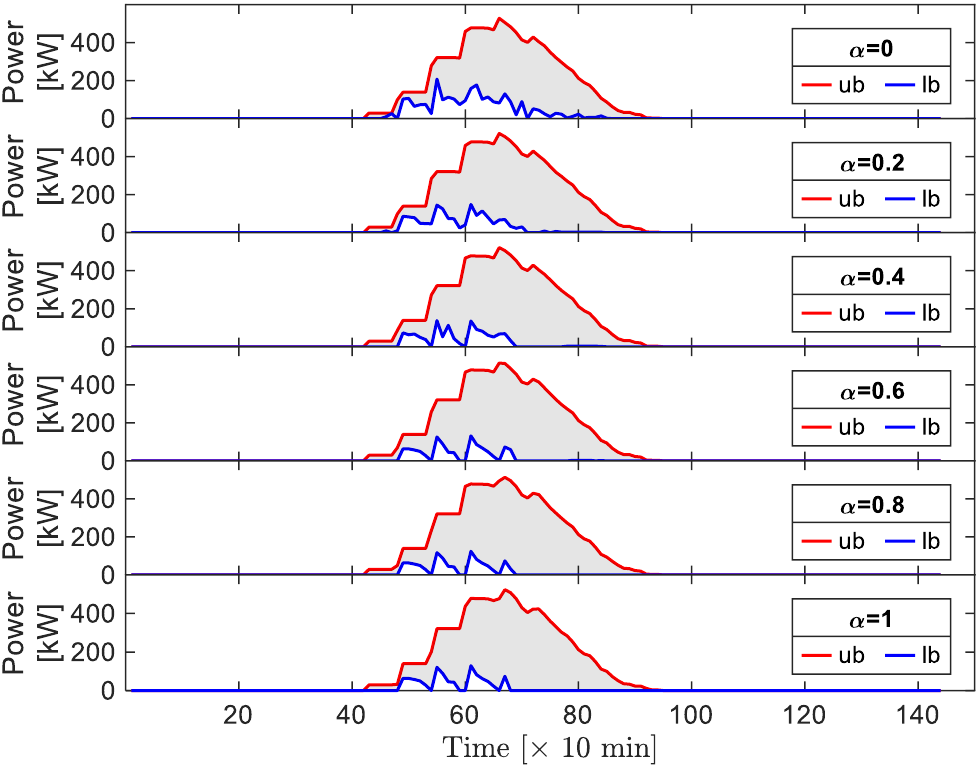} \\
  \caption{The aggregate EV power flexibility under different $\alpha$.}
  \label{fig:flex}
\end{figure}

The above result is obtained under the random dispatch ratio $\alpha_t$ (see Fig. \ref{fig:alpha}).
In fact, the dispatch ratio can affect the actual charging power of each EV and further affect their aggregate power flexibility. Therefore, it is interesting to investigate the impact of $\alpha_t$ on the aggregate EV power flexibility. Here, we use a uniform $\alpha$ over time, i.e., $\alpha_t=\alpha, \forall t$. We change $\alpha$ from 0 to 1 and record the total power flexibility value in Fig.~\ref{fig:diff_alpha}. As seen, the total power flexibility value depends on the dispatch ratio $\alpha$.
Generally, a larger $\alpha$ leads to a larger power flexibility value. However, this increase is nonlinear. When the dispatch ratio $\alpha$ exceeds 0.9, the total power flexibility value no longer increases. In the extreme case when $\alpha=0$, the total power flexibility value is 262.3 USD, which is still larger than that of the greedy algorithm B1. In addition, it can be concluded that if the average dispatch ratio $\alpha$ is greater than 0.15, then the proposed online algorithm more likely outperforms the offline method. This demonstrates the advantage of the proposed algorithm.
We also present the aggregate EV power flexibility region under different $\alpha$ in Fig. \ref{fig:flex}. As $\alpha$ decreases, the aggregate EV power flexibility region gradually narrows. This is because, under a low dispatch ratio, the EVs are charged at a low charging rate and are more likely to fail to meet the charging requirements; hence, the lower bound of the aggregate EV power flexibility region is raised to ensure that the EVs can meet the charging requirements in the remaining time slots.

\subsection{Impact of Parameters}
\label{secV-D}
\begin{figure}[!htbp]
  \centering
  \includegraphics[width=0.45\textwidth]{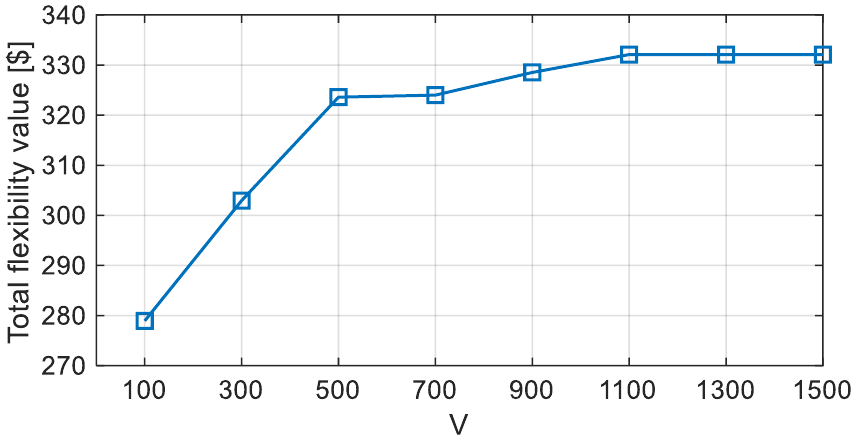} \\
  \caption{The impact of $V$ on the total power flexibility value.}
  \label{fig:diffv}
\end{figure}

According to \eqref{equ:driftPlusP}, the parameter $V$ controls the trade-off between stabilizing the queues and maximizing the total power flexibility value. Here, we change the value of $V$ to investigate its impact on the total power flexibility value. As shown in Fig.~\ref{fig:diffv}, the total power flexibility value becomes larger with an increasing $V$.

\begin{figure}[!htbp]
  \centering
  \includegraphics[width=0.49\textwidth]{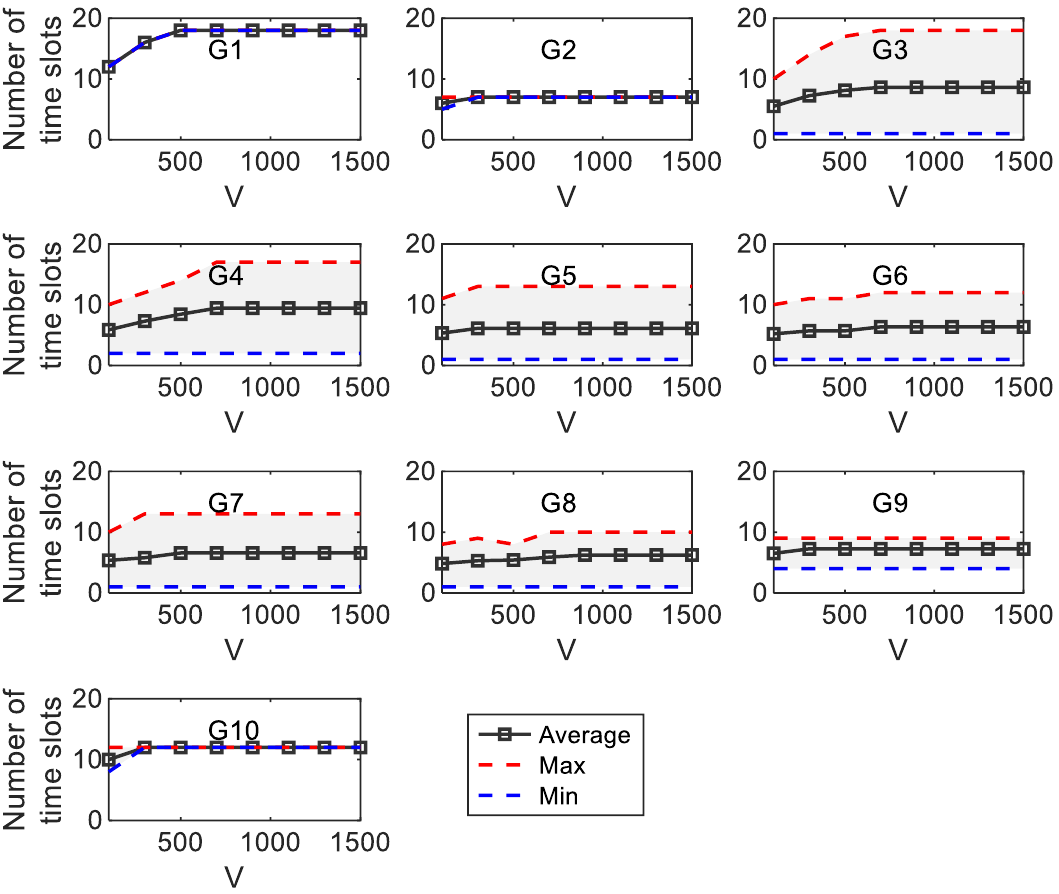} \\
  \caption{The impact of $V$ on the number of time slots taken to reach the charging requirement.}
  \label{fig:delay}
\end{figure}

Fig. \ref{fig:delay} depicts the impact of $V$ on the number of time slots needed for EVs to meet their charging requirement (SOC $=0.5$). We calculate the maximum/minimum/average number of time slots needed for the EVs in each group. As seen, with the growth of $V$, the number of time slots needed slightly increases. This is because a larger $V$ means putting more emphasis on maximizing the total power flexibility, which results in a reduced lower bound $\check{x}_{g,t}$ of the aggregate EV power flexibility region. Consequently, the charging time needed to reach the required SoC becomes longer. 

\begin{figure}[!htbp]
  \centering
  \includegraphics[width=0.47\textwidth]{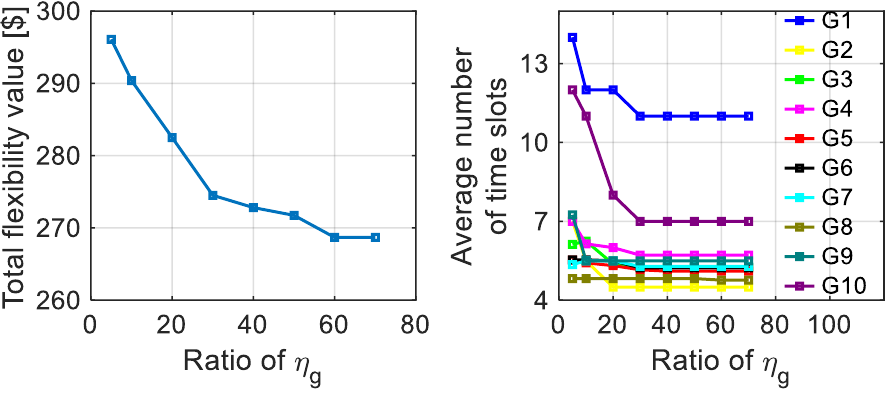} \\
  \caption{The impact of $\eta_g$ on the total flexibility value and the average number of time slots taken to reach the charging requirement.}
  \label{fig:eta}
\end{figure}

Fig. \ref{fig:eta} shows the impact of $\eta_g$ on the total flexibility value and the number of time slots needed for EVs to meet their required charging SOC $=0.5$. We can find that a larger $\eta_g$ results in a lower total power flexibility value and less number of time slots taken. This is because a larger $\eta_g$ forces the virtual delay-aware queue $\check{Z}_g$ to grow more rapidly, allowing EVs to get charged quickly. Meanwhile, the power flexibility is sacrificed.

\subsection{Sensitivity and Scalability Analysis}
\label{secV-E}
\begin{figure}[!htbp]
  \centering
  \includegraphics[width=0.3\textwidth]{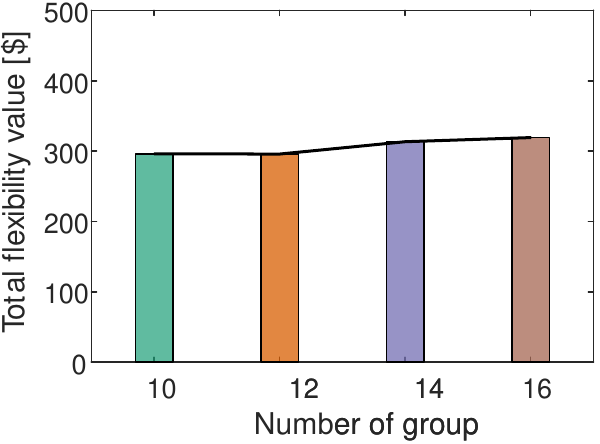} \\
  \caption{Impact of group numbers on the total power flexibility value.}\label{fig:NRGvsfv}
\end{figure}

\begin{figure}[!htbp]
  \centering
  \includegraphics[width=0.45\textwidth]{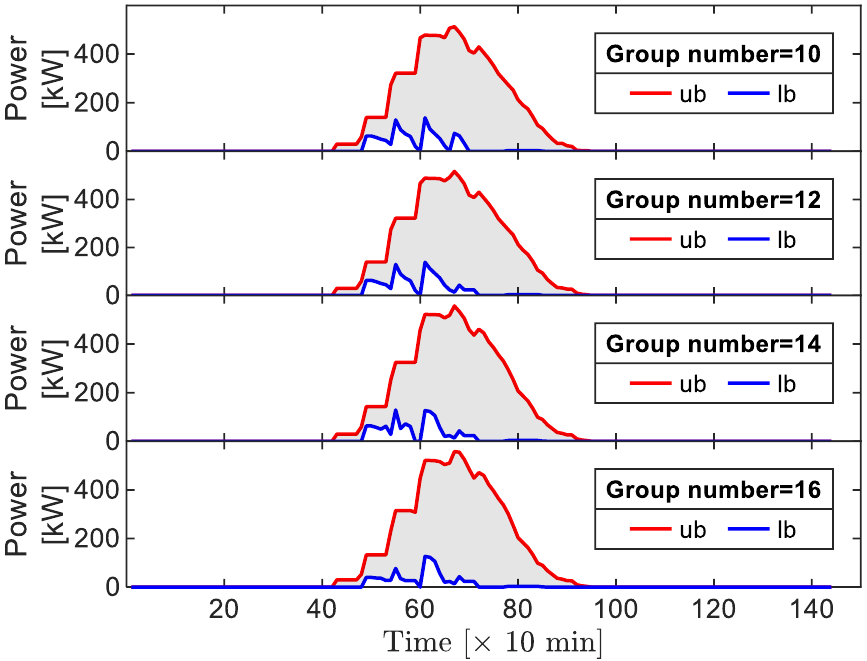} \\
  \caption{Aggregate power flexibility region under different number of groups.}\label{fig:NRGvsagg}
\end{figure}

In the following, we conduct sensitivity and scalability analysis. First, the impact of the number of groups on the total flexibility value is investigated. Under different cases, the total number of EVs remains the same. We split the EVs that are originally in one group (G3, G4, G5, G6, G7, G8) equally into two groups and modify their original charging delay value $R_g$ a little bit into $R_g-1$ and $R_g+1$. Fig. \ref{fig:NRGvsfv} shows that the total power flexibility value increases slightly as the number of groups increases. Fig. \ref{fig:NRGvsagg} gives the aggregate power flexibility regions for the cases with 10, 12, 14, and 16 groups, respectively. It can be seen that a smaller group number leads to a smaller aggregate power flexibility region and total power flexibility value. This is because a larger number of groups results in a smaller $\hat Q_{g,t}$, so the objective function of \textbf{P3} emphasizes more on the flexibility value maximization.

\begin{figure}[!htbp]
  \centering
  \includegraphics[width=0.45\textwidth]{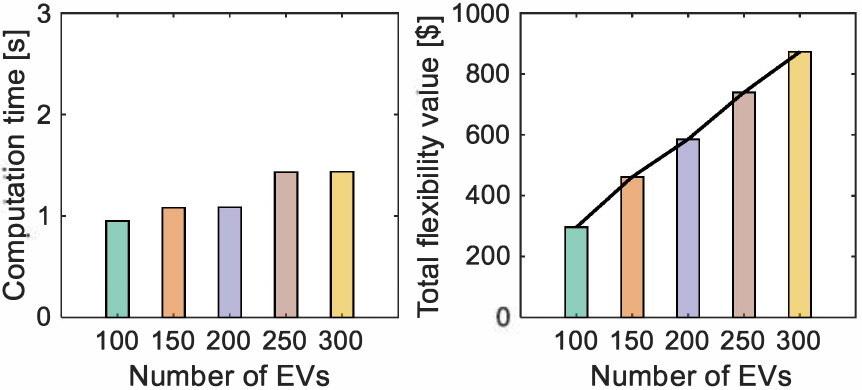} \\
  \caption{Computation time and total power flexibility value under different number of EVs.}
  \label{fig:computeTime}
\end{figure}

To show the scalability of the proposed online algorithm, we test its performance under different numbers of EVs, ranging from 100 to 300. The left-hand side of Fig. \ref{fig:computeTime} shows the change of the computation time. For the case with 100 EVs, the algorithm only takes 0.95s. In addition, we can find that the computation time does not significantly increase as the number of EVs grows, demonstrating the scalability of the proposed online algorithm. The right-hand side of Fig. \ref{fig:computeTime} shows the total power flexibility value under different numbers of EV. With more EVs, the total power flexibility value almost increases proportionally. This is consistent with the intuition that the more EVs, the higher the flexibility.

\subsection{More Challenging Situations}
\label{secV-F}
\begin{figure}[!htbp]
  \centering
  \includegraphics[width=0.45\textwidth]{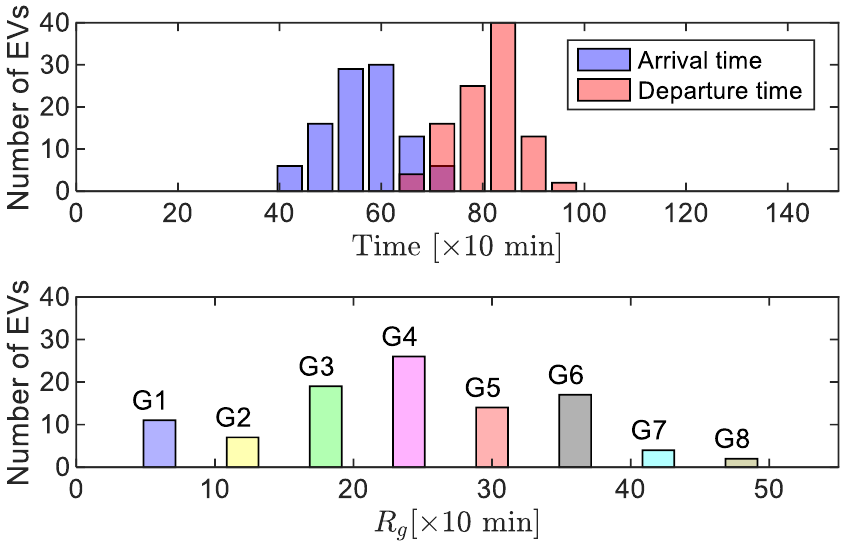}\\
  \caption{Distribution of EV arrival time and departure time, and distribution of $R_g$.}\label{fig:time_case2}
\end{figure}

We further test the performance of the proposed algorithm under more challenging situations. We consider a higher desired SOC value of 0.7 and shift the departure time earlier, i.e., from the original distribution $\mathcal{N}(\text{18:00}, (1.2hr)^2)$ to the distribution $\mathcal{N}(\text{14:00},(1.2hr)^2)$. The allowed charging delays $R_g,\forall g$ are shown in Fig. \ref{fig:time_case2}, which are considerably shorter. As seen, 63\% of the EVs (those in G1-G4) stay for 4 hours or less, with EVs in G1 having the shortest stay time of only 1 hour. In addition, we randomly set the EV initial SOC from a normal distribution $\mathcal{N}(0.4, 0.01)$ \cite{liu2018dispatch}. But it is worth noting that the assumption of a uniform distribution for the initial SOC of EVs is also widely used in the EV charging literature such as \cite{jin2014optimized,wei2018charging,koufakis2020offline}.

\begin{figure}[!htbp]
  \centering
  \includegraphics[width=0.5\textwidth]{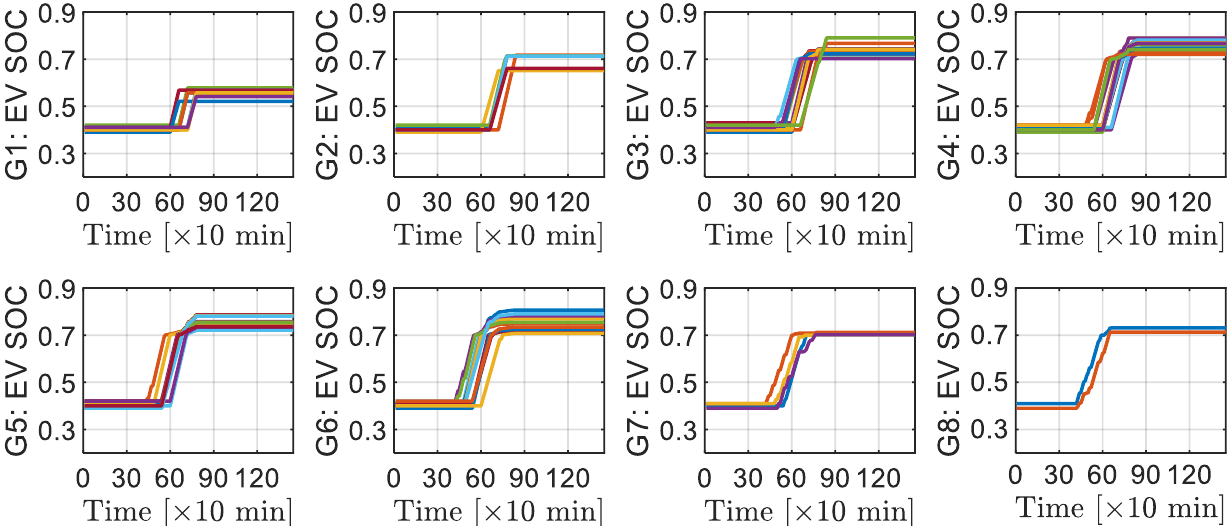}\\
  \caption{EV SOC of different groups.}\label{fig:evsoc_case2}
\end{figure}

\begin{figure}[!htbp]
  \centering
  \includegraphics[width=0.5\textwidth]{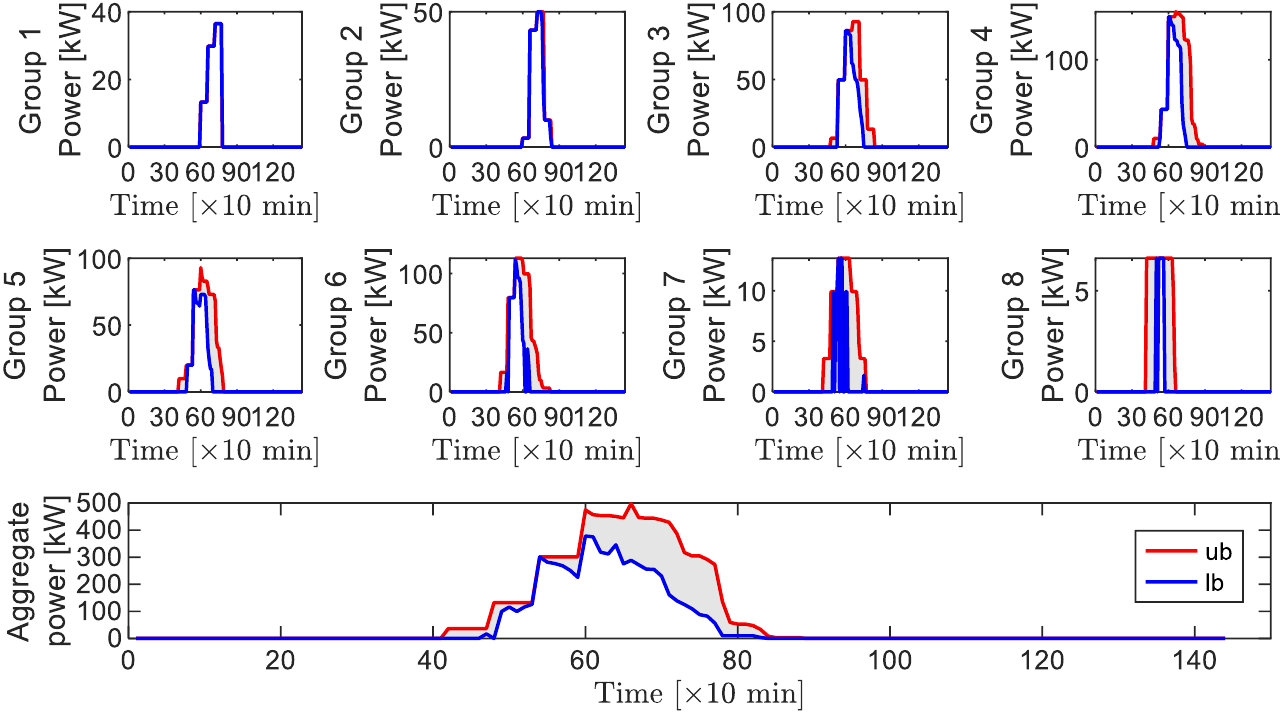}\\
  \caption{Aggregate EV power flexibility region.}\label{fig:flex_prop_case2}
\end{figure}

Fig. \ref{fig:evsoc_case2} shows the EV charging SOC curves under the more challenging setting above. As seen, all EVs in G1 and most EVs in G2 cannot reach the desired SOC of 0.7. This is because their minimal required charging time (calculated by the charging-as-soon-as-possible method) is larger than their allowed charging delay $R_g$. Therefore, even under the offline method, those EVs cannot be fully charged. As shown in Fig. \ref{fig:flex_prop_case2}, G1 cannot provide power flexibility and G2 offers very little. Meanwhile, the aggregate power flexibility region (Fig. \ref{fig:flex_prop_case2}, bottom) is much narrower than that in Fig. \ref{fig:flex_prop}. This demonstrates that a lower desired SOC and shorter charging time leads to a narrower power flexibility region, and vice versa. In addition, the total flexibility value obtained by the proposed method is 115.4 USD which is still higher than that by the offline method (113.6 USD).

\section{Conclusion}\label{sec:conclu}
With the proliferation of EVs, it is necessary to better utilize their charging power flexibility, making them valuable resources rather than burdens on the power grid. This paper proposes a real-time feedback based online aggregate EV power flexibility characterization method. It can output the aggregate EV flexibility region for each time slot in an online manner, with the total power flexibility over time similar to that of the offline counterpart. We prove that by choosing an aggregate dispatch strategy within the obtained flexibility region for each time slot, the corresponding disaggregated EV control strategies allow all EVs to satisfy their charging requirements. The case studies demonstrate the effectiveness and benefits of the proposed method. We have the following findings:

1) The proposed online algorithm can even outperform the offline method in some cases since it can utilize the up-to-date dispatch information via real-time feedback.

2) The operator's dispatch decision can influence the future aggregate EV power flexibility. Particularly, a higher dispatch ratio results in larger flexibility regions in future time slots.


3) A larger parameter $V$ leads to a higher total flexibility but a longer time to meet the EV charging requirements.

Future research may take into account the conflicting interests between the operator, aggregators, and EVs when deriving the aggregate flexibility region.
In addition, how to fully protect the privacy of EVs is another possible direction.

\appendices
\makeatletter
\@addtoreset{equation}{section}
\@addtoreset{theorem}{section}
\makeatother
\setcounter{equation}{0}  
\renewcommand{\theequation}{A.\arabic{equation}}
\renewcommand{\thetheorem}{A.\arabic{theorem}}

\section{Proof of Proposition \ref{prop-1}}
\label{appendix-1}
Let $\{p_{d,t},\forall t\}$ be the aggregate power trajectory. For each time slot $t \in \mathcal{T}$, since $p_{d,t} \in [\check p_{d,t}^*, \hat p_{d,t}^*]$, we can define an auxiliary coefficient:
\begin{align}
    \beta_t:= \frac{\hat p_{d,t}^*-p_{d,t}}{\hat p_{d,t}^*-\check p_{d,t}^*} \in [0,1]
\end{align}
so that $p_{d,t}=\beta_t \check{p}^*_{d,t}+(1-\beta_t)\hat{p}^*_{d,t}$. Then, we can construct a feasible EV dispatch strategy by letting
\begin{subequations}
\begin{align}
    p^c_{v,t}&=\beta_t\check{p}^{c*}_{v,t}+(1-\beta_t)\hat{p}^{c*}_{v,t},\\
    e_{v,t}&=\beta_t\check{e}^{c*}_{v,t}+(1-\beta_t)\hat{e}^{c*}_{v,t}.
\end{align}
\end{subequations}
for all time slots $t \in \mathcal{T}$.

We prove that it is a feasible EV dispatch strategy as follows, 
\begin{align}
    p_{d,t}=~ & \beta_t\check p_{d,t}^* + (1-\beta_t) \hat p_{d,t}^* \nonumber\\
    =~ &  \beta_t \sum_{v \in \mathcal{V}} \check p_{v,t}^{c*} + (1-\beta_t) \sum_{v \in \mathcal{V}} \hat p_{v,t}^{c*} \nonumber\\
    = ~ & \sum_{v \in \mathcal{V}} \left[\beta_t \check p_{v,t}^{c*}+(1-\beta_t) \hat p_{v,t}^{c*}\right] \nonumber\\
    = ~ & \sum_{v \in \mathcal{V}} p^c_{v,t}
\end{align}
Hence, constraint \eqref{equ:ubpdi} holds for $p_{d,t}$ and $p^c_{v,t},\forall v$. Similarly, we can prove that constraints \eqref{equ:ubpd}-\eqref{equ:uberanges} and \eqref{eq:evava} are met. Therefore, we have constructed a feasible EV dispatch strategy, which completes the proof. \hfill$\blacksquare$

\setcounter{equation}{0}  
\renewcommand{\theequation}{C.\arabic{equation}}
\renewcommand{\thetheorem}{C.\arabic{theorem}}

\setcounter{equation}{0}  
\renewcommand{\theequation}{B.\arabic{equation}}
\renewcommand{\thetheorem}{B.\arabic{theorem}}

\section{Proof of Proposition \ref{prop-2}}
\label{appendix-2}
Here, we use the contradiction. If a charging request $\hat{a}_{g,t}$ arrives in time slot $t$ cannot be fulfilled on or before time slot $t+\hat{\delta}_{g,max}$. Then, queue $\hat{Q}_{g,\tau}>0$ always holds for $\tau\in[t+1,...,t+\hat{\delta}_{g,max}]$. Thus, we have $\mathbb{I}_{\hat{Q}_{g,\tau}>0}=1,\forall \tau$. According to delay virtual queue dynamics \eqref{eq:zgt-ub}, for all $\tau\in[t+1,...,t+\hat{\delta}_{g,max}]$, we have
\begin{gather}
    \hat{Z}_{g,\tau+1} \geq \hat{Z}_{g,\tau}+\frac{\eta_g}{R_g}-\hat{x}_{g,\tau}, \forall g, \forall t.
\end{gather}
By summing the above inequalities over $\tau\in[t+1,...,t+\hat{\delta}_{g,max}]$, we have
\begin{gather}
    \hat{Z}_{g,t+\hat{\delta}_{g,max}+1}-\hat{Z}_{g,t+1} \geq    \frac{\eta_g}{R_g}\hat{\delta}_{g,max}+\sum_{\tau=t+1}^{t+\hat{\delta}_{g,max}}(-\hat{x}_{g,\tau}).
\end{gather}
Since $\hat{Z}_{g,t+\hat{\delta}_{g,max}+1}\leq \hat{Z}_{g,max}$ and $\hat{Z}_{g,t+1}\geq 0$, we have
\begin{gather}
    \hat{Z}_{g,max} \geq \frac{\eta_g}{R_g}\hat{\delta}_{g,max}+\sum_{\tau=t+1}^{t+\hat{\delta}_{g,max}}(-\hat{x}_{g,\tau}).
\end{gather}
Since the charging tasks are processed in a first-in-first-out manner, and the charging request is not fulfilled by $t+\hat{\delta}_{g,max}$, we have
\begin{equation}
    \sum_{\tau=t+1}^{t+\hat{\delta}_{g,max}}(\hat{x}_{g,\tau})<\hat{Q}_{g,max}
\end{equation}
Combining the above two inequalities, we obtain
\begin{gather}
    \hat{Z}_{g,max} > \frac{\eta_g}{R_g}\hat{\delta}_{g,max}-\hat{Q}_{g,max},
\end{gather}
which implies
\begin{equation}
    \hat{\delta}_{g,max}<\frac{(\hat{Q}_{g,max}+\hat{Z}_{g,max})R_{g}}{\eta_g}.\\
\end{equation}
However, this result contradicts the definition of $\hat{\delta}_{g,max}$ in \eqref{eq:dgmax-ub}. Therefore, the worst case delay should be less than or equal to $\hat{\delta}_{g,max}$ as defined in \eqref{eq:dgmax-ub}.

The proof of \eqref{eq:dgmax-lb} follows a similar procedure, and we omit it here for brevity. $\hfill \blacksquare$

\setcounter{equation}{0}  
\renewcommand{\theequation}{C.\arabic{equation}}
\renewcommand{\thetheorem}{C.\arabic{theorem}}

\section{Proof of Proposition \ref{prop-3}}
\label{appendix-3}
Denote the optimal solution of \textbf{P3} as $\hat{x}_{g,t}^{*}$ and $\check{x}_{g,t}^{*}$, and the optimal solution of \textbf{P2} by $\hat{x}_{g,t}^{off}$ and $\check{x}_{g,t}^{off}$. According to \eqref{equ:dppInequ}, we have
\begin{equation}\label{eq:gap}
\begin{aligned}
&\mathbb{E}\left[\Delta(\boldsymbol{\Theta}_t)|\boldsymbol{\Theta}_t\right]+V\mathbb{E}[-F_t^{*}|\boldsymbol{\Theta}_t]\\
&\leq A+V\mathbb{E}[-F_t^{*}|\boldsymbol{\Theta}_t]+\sum\limits_{g\in{\mathcal{G}}}\hat{Q}_{g,t}\mathbb{E}\left[\hat{a}_{g,t}-\hat{x}_{g,t}^{*}|\boldsymbol{\Theta}_t\right]\\
&+\sum\limits_{g\in{\mathcal{G}}}\check{Q}_{g,t}\mathbb{E}\left[\check{a}_{g,t}-\check{x}_{g,t}^{*}|\boldsymbol{\Theta}_t\right]+\sum\limits_{g\in{\mathcal{G}}}\hat{Z}_{g,t}\mathbb{E}\left[-\hat{x}_{g,t}^{*}|\boldsymbol{\Theta}_t\right]\\
&+\sum\limits_{g\in{\mathcal{G}}}\check{Z}_{g,t}\mathbb{E}\left[-\check{x}_{g,t}^{*}|\boldsymbol{\Theta}_t\right],\\
&\leq A+V\mathbb{E}[-F_t^{off}|\boldsymbol{\Theta}_t]+\sum\limits_{g\in{\mathcal{G}}}\hat{Q}_{g,t}\mathbb{E}\left[\hat{a}_{g,t}-\hat{x}_{g,t}^{off}|\boldsymbol{\Theta}_t\right]\\
&+\sum\limits_{g\in{\mathcal{G}}}\check{Q}_{g,t}\mathbb{E}\left[\check{a}_{g,t}-\check{x}_{g,t}^{off}|\boldsymbol{\Theta}_t\right]+\sum\limits_{g\in{\mathcal{G}}}\hat{Z}_{g,t}\mathbb{E}\left[-\hat{x}_{g,t}^{off}|\boldsymbol{\Theta}_t\right]\\
&+\sum\limits_{g\in{\mathcal{G}}}\check{Z}_{g,t}\mathbb{E}\left[-\check{x}_{g,t}^{off}|\boldsymbol{\Theta}_t\right].
\end{aligned}
\end{equation}
By summing the above inequality \eqref{eq:gap} over time slots $t\in\{1,2,\ldots,T\}$, dividing both sides by $VT$, and letting $T \to \infty$, we have
\begin{align}
   & \lim_{T \to \infty} \frac{1}{T}\left(\mathbb{E}[L(\boldsymbol{\Theta}_{T+1})]-\mathbb{E}[L(\boldsymbol{\Theta}_1)]\right) + \lim_{T \to \infty} \frac{V}{T} \sum_{t=1}^T \mathbb{E}(-F_t^{*}) \nonumber\\
  & \le  A + \lim_{T \to \infty} \frac{V}{T} \sum_{t=1}^T \mathbb{E}(-F_t^{off}).
\end{align}

The is based on the fact that
\begin{align}
    \lim\limits_{T\rightarrow\infty}\frac{1}{T}\sum\limits_{t=1}^{T}\mathbb{E}\left[\hat{a}_{g,t}-\hat{x}_{g,t}^{off}|\boldsymbol{\Theta}_t\right]\leq 0,\\
    \lim\limits_{T\rightarrow\infty}\frac{1}{T}\sum\limits_{t=1}^{T}\mathbb{E}\left[\check{a}_{g,t}-\check{x}_{g,t}^{off}|\boldsymbol{\Theta}_t\right]\leq 0,\\
    \lim\limits_{T\rightarrow\infty}\frac{1}{T}\sum\limits_{t=1}^{T}\mathbb{E}\left[-\hat{x}_{g,t}^{off}|\boldsymbol{\Theta}_t\right]\leq 0,\\
    \lim\limits_{T\rightarrow\infty}\frac{1}{T}\sum\limits_{t=1}^{T}\mathbb{E}\left[-\check{x}_{g,t}^{off}|\boldsymbol{\Theta}_t\right]\leq 0,
\end{align}
which is due to constraints \eqref{eq:Qgub-lim}-\eqref{eq:xglb}.

Since $L(\boldsymbol{\Theta}_{T+1})$ and $L(\boldsymbol{\Theta}_1)$ are finite, we have
\begin{align*}
\underbrace{\lim\limits_{T\rightarrow\infty}\frac{1}{T}\sum\limits_{t=1}^{T}\mathbb{E}[F_t^{*}]}_{\mathcal{F}^{*}}\ge -\frac{A}{V}+\underbrace{\lim\limits_{T\rightarrow\infty}\frac{1}{T}\sum\limits_{t=1}^{T}\mathbb{E}[F_t^{off}]}_{\mathcal{F}^{off}}.
\end{align*}
Moreover, it is easy to know $\mathcal{F}^{off} \ge \mathcal{F}^{*}$. $\hfill \blacksquare$

\bibliographystyle{IEEEtran}
\bibliography{PaperRef}

\end{document}